\documentclass[11pt]{amsart}
\usepackage{amssymb}
\usepackage{latexsym}
\usepackage{amsmath}
\usepackage{amsthm}
\usepackage{amsfonts}
\usepackage{color}
\usepackage{pdfsync}
\usepackage[usenames,dvipsnames]{pstricks}
\usepackage{epsfig}
\usepackage{pst-grad} 
\usepackage{pst-plot} 

\newcommand{\R}{\mathbb{R}}

\theoremstyle{plain}
\newtheorem{defi}{Definition}[section]
\newtheorem{prop}[defi]{Proposition}
\newtheorem{teo}[defi]{Theorem}
\newtheorem{cor}[defi]{Corollary}
\newtheorem{lema}[defi]{Lemma}

\newtheorem{remark}[defi]{Remark}

\theoremstyle{definition}

\theoremstyle{remark}

\numberwithin{equation}{section}

\begin{document}

\title[]{Existence and Uniqueness for Integro-Differential Equations with Dominating Drift Terms.}

%

\author[]{Erwin Topp}
\address{
Erwin Topp - 
Departamento de Ingenier\'\i a Matem\'atica (UMI 2807 CNRS), Universidad de Chile, Casilla 170, Correo 3, Santiago, Chile (etopp@dim.uchile.cl).}

\address{and}

\address{Laboratoire de Math\'ematiques et Physique Th\'eorique (CNRS UMR 6083), F\'ederation Denis Poisson, 
Universit\'e Fran\c{c}ois Rabelais, Parc de Grandmont, 37200, Tours, France (Erwin.Topp@lmpt.univ-tours.fr).}



\date{\today}

\begin{abstract} 
In this paper we are interested on the well-posedness of Dirichlet problems associated to integro-differential 
elliptic operators of order $\alpha < 1$ in a bounded smooth domain $\Omega$ . The main difficulty arises because of losses of the boundary condition 
for sub and supersolutions due to the lower diffusive effect of the elliptic operator compared with the drift term. We consider the 
notion of viscosity solution with generalized boundary conditions, concluding strong comparison principles in $\bar{\Omega}$ 
under rather general assumptions over the drift term. As a consequence, existence and uniqueness of solutions in 
$C(\bar{\Omega})$ is obtained via Perron's method.
\end{abstract}

\maketitle

\section{Introduction.}

Let $\Omega \subset \R^n$ be a bounded domain with $C^2$ boundary, $\varphi: \Omega^c \to \R$ a bounded continuous function, 
$H: \bar{\Omega} \times \R^n \to \R$ a continuous function, $\lambda \geq 0$
and $\alpha \in (0,1)$. The main interest of this paper is the study of the Dirichlet problem
\begin{eqnarray}
\label{eq} \lambda u - \mathcal{I}[u] + H(x, Du) \ = & 0 & \quad \mbox{in} \ \Omega, \\
\label{Dirichlet} u \ = & \varphi & \quad \mbox{on} \ \Omega^c.
\end{eqnarray}

The term $\mathcal{I}$ represents an integro-differential operator defined as follows: Consider $\Lambda > 0$ and
$K \in L^\infty(\R^n)$ a nonnegative function such that $||K||_{L^\infty(\R^n)} \leq \Lambda$. For 
$x \in \R^n$ and $\phi: \R^n \to \R$ bounded and sufficiently regular at $x$, we define
\begin{eqnarray*}\label{operator}
\mathcal{I}[\phi](x) = \int \limits_{\R^n} [\phi(x + z) - \phi(x)] K^\alpha(z) dz,
\end{eqnarray*} 
where $K^\alpha(z) := K(z)|z|^{-(n + \alpha)}$ for $z \neq 0$.  
A particular case
is when $K$ is equal to a well-known constant $C_{n,\alpha} > 0$ since, in that case, $-\mathcal{I} = (-\Delta)^{\alpha/2}$, 
the fractional Laplacian of order $\alpha$ (see~\cite{Hitch}). 

Note that condition~\eqref{Dirichlet} imposes the unknown function $u$ is equal to $\varphi$ in $\Omega^c$. We incorporate 
this exterior data into the equation through the nonlocal term and therefore the problem can be as well reformulated as : finding a function 
$u \in C(\bar{\Omega})$ satisfying, for all $x \in \Omega$, the following equality in the viscosity sense 
\begin{equation}\label{eqintro}
\begin{split}
\lambda u(x) & - \int \limits_{x + z \in \bar{\Omega}} [u(x + z) - u(x)] K^\alpha(z)dz \\
& - \int \limits_{x + z \notin \bar{\Omega}} [\varphi(x + z) - u(x)] K^\alpha(z)dz + H(x, Du(x)) = 0,
\end{split}
\end{equation}
together with an appropriate notion of boundary condition to be precised later on.
This perspective of the problem is a way to deal with the eventual existence of discontinuities at the boundary, namely $u$ may be different 
to $\varphi$ at some points of $\partial \Omega$. 

Related to this, we cite 
the work of Barles, Chasseigne and Imbert~\cite{Barles-Chasseigne-Imbert}, where the authors study the Dirichlet 
problem for a large variety of integro-differential equations, proving in particular that, in the absence of second-order 
terms, there is no loss of the boundary condition provided the order of the fractional operator is $\alpha \geq 1$, no matter the behavior 
of the drift terms is. On the contrary, loss of the boundary condition may occur in problems like~\eqref{eq}-\eqref{Dirichlet} when $\alpha < 1$ 
and the first-order term has certain geometric disposition at the boundary, providing an explicit example of this phenomena.
Therefore, our aim is to get the well-posedness of problem~\eqref{eq}-\eqref{Dirichlet} in this context 
of loss of the boundary condition, arising under this competitive relation between a ``weak'' diffusive effect of the nonlocal term $\mathcal{I}$ and 
the first-order transport effect of $H$. Thus, the knowledge of the nature of the first-order term is determinant  in this task.

For this, we assume $H$ has a \textsl{Bellman form} defined in the following way: Let $\mathcal{B}$ be a metric compact space, 
$b: \bar{\Omega} \times \mathcal{B} \to \R^n$, $f: \bar{\Omega} \times \mathcal{B} \to \R$ continuous and bounded. 
For $x \in \bar{\Omega}$ and $p \in \R^n$, we denote
\begin{eqnarray}\label{Hamiltonian}
H(x, p) = \sup \limits_{\beta \in \mathcal{B}} \{ -b_\beta(x) \cdot p - f_\beta(x)\}
\end{eqnarray}
where we adopt the notation $b_\beta(x) = b(x, \beta)$ for $(x, \beta) \in \bar{\Omega} \times \mathcal{B}$ and 
in the same way for $f$. Hamiltonians as~\eqref{Hamiltonian} arise in the study of Hamilton-Jacobi equations associated to
optimal exit time problems (see~\cite{Bardi-Capuzzo},~\cite{Barles-book},~\cite{Fleming-Soner}). In such problems, the study of $H$ is made 
through $b$, concluding that its behavior at the boundary may affect in a sensible way the properties of the equation.

An illustrative model which is similar to ours is the exhaustively studied \textsl{degenerate second-order elliptic} case 
(see~\cite{Keldysh} and the monographs~\cite{Radkevich1},~\cite{Radkevich2} for a complete survey of the also called 
\textsl{equations with nonnegative characteristic form}). For a
continuous function $a$ with values in the set of nonnegative matrices, $f$ and $\varphi$ continuous real valued functions, $\lambda \geq 0$ 
a constant and $b$ a continuous vector field, a linear model equation
\begin{eqnarray}\label{linearsecondorder}
\left \{ \begin{array}{rcll} - {\rm Tr}(a(x) D^2u) - b(x) \cdot Du + \lambda u &=& f(x) & \mbox{in } \ \Omega \\
u &=& \varphi & \mbox{on} \ \partial \Omega, \end{array} \right .
\end{eqnarray}
is said to be degenerate elliptic if the matrix $a(x)$ has null eigenvalues for some $x \in \bar{\Omega}$, being of special interest the case 
when the degeneracy points are on the boundary. It is known (see~\cite{Freidlin}) that if a solution $u$ of~\eqref{linearsecondorder} is such 
that $u(x_0) > \varphi(x_0)$ for some $x_0 \in \partial \Omega$, then necessarily
\begin{eqnarray}\label{nonregularpoint}
\left \{ \begin{array}{l} a(x_0) Dd(x_0) \cdot Dd(x_0) = 0, \ and\\ -{\rm Tr}(a(x_0) D^2d(x_0)) - b(x_0) \cdot Dd(x_0) \leq 0. \end{array} \right .
\end{eqnarray}

Since $Dd(x_0)$ agrees with the inward unit normal of $\partial \Omega$
at $x_0$, first condition establishes the degeneracy of the second-order operator, depicting the absence of diffusion in the normal direction at 
$x_0$. Second condition shows that certain geometric disposition of the drift term at the boundary is necessary to loss the boundary condition.
Even under this difficulty, comparison principle for bounded sub and supersolutions can be obtained for degenerate elliptic
second-order problems. For instance, in the context of viscosity solutions with \textsl{generalized boundary conditions}, 
Barles and Burdeau in~\cite{Barles-Burdeau} prove a comparison result for a variety of 
quasilinear problems associated to Bellman-type Hamiltonians as~\eqref{Hamiltonian}. Related results can be also found 
in~\cite{Barles-DaLio},~\cite{Barles-Rouy},~\cite{DaLio}.

Condition~\eqref{nonregularpoint} can be understood using the approach of stochastic exit time problems.
Assume $\sigma$ is a continuous function 
such that $a = 1/2 \sigma \sigma^t$, denote $(W_t)_t$ the standard Brownian motion, let $x \in \bar{\Omega}$ and consider the stochastic 
differential equation
\begin{equation}\label{X}
dX_t^x = b(X_t^x)dt + \sigma(X_t^x)dW_t, \quad X_0^x = x.  
\end{equation}

If we denote $\tau$ the exit time from $\Omega$ of $(X_t)_t$ and $\mathbb{E}_x$ the conditional expectation with respect to 
the event $\{ X_0 = x \}$, it is known (see~\cite{Freidlin} for instance) the \textsl{value function}
\begin{equation}\label{payoff}
\mathbb{U}(x) = \mathbb{E}_x \Big{\{} \int \limits_{0}^{\tau} f(X_s)e^{-\lambda s}ds + \varphi(X_\tau) e^{-\lambda \tau} \Big{\}}
\end{equation}
is a solution to~\eqref{linearsecondorder}. Thus, for a point $x_0$ close to the boundary, the degeneracy of the diffusion driven by $\sigma$ at 
$x_0$ allows $b$ ``to push'' the trajectories inside the domain and then $\mathbb{U}(x_0)$ is not necessarily equal to $\varphi(x_0)$,
explaining the relation among~\eqref{nonregularpoint} and the loss of the boundary condition.

At this point, we introduce the key assumption
\begin{equation}\tag{{\bf E}}\label{ellipticity}
\exists \ c_1, c_2 > 0 \quad \mbox{such that} \quad c_1 \leq K(z) \quad \mbox{for any} \ |z| \leq c_2,
\end{equation} 
which has to be understood as an \textsl{uniform ellipticity condition} over $\mathcal{I}$ in the sense of Caffarelli and 
Silvestre~\cite{Caffarelli-Silvestre1}. By~\eqref{ellipticity}, $\mathcal{I}$ is the infinitesimal generator 
of a L\'evy process $(Z_t)_t$ whose a.a. paths have infinite number of jumps in each interval of time, but finite variation by the assumption
$\alpha < 1$ (see~\cite{Sato}). Then, if $(\beta_t)_t$ is a stochastic process $(Z_t)_t-$adapted with values in $\mathcal{B}$ and 
$x \in \bar{\Omega}$, the SDE
\begin{equation}\label{trajectories}
dX_t^\beta = b(X_t^\beta, \beta_t)dt + dZ_t, \quad X_0^\beta = x,
\end{equation}
allows us to give an stochastic interpretation to our integro-differential problem analogous to the second-order case, playing
the role of~\eqref{X}. 
We recall that losses of boundary conditions are present in both approaches. However, 
we would like to emphasize that the reasons for these phenomenas are qualitatively different. 
In the second-order case, they only arise when the diffusion \textsl{is degenerate}, preventing the trajectories driven by~\eqref{X} 
reach the boundary. On the contrary, the \textsl{uniform ellipticity condition}~\eqref{ellipticity} allows 
the trajectories defined in~\eqref{trajectories} make small jumps in all directions, in particular, they can jump outside $\Omega$
when they start close to the boundary. This last feature intuitively should insure the boundary condition, but the transport effect of the drift term
is strong enough to produce loss of boundary conditions, mainly due to the assumption $\alpha < 1$.

The above discussion gives rise to the main purpose of this paper, which is to establish a viscosity comparison principle for 
bounded sub and supersolutions of a class of integro-differential problems which seem to be uniformly elliptic, but where the functions can be 
discontinuous at the boundary.

We finish this introduction with a useful interpretation of our problem. For $x \in \bar{\Omega}$ and $u$ a suitable function, define the
\textsl{censored nonlocal operator} as
\begin{equation}\label{censoredI}
\mathcal{I}_\Omega[u](x) := \int \limits_{x + z \in \bar{\Omega}} [u(x + z) - u(x)] K^\alpha(z)dz.
\end{equation}

For $x \in \Omega$, we also define
\begin{equation*}\label{unboundedcoef}
\bar{\lambda}(x) = \lambda + \int \limits_{x + z \notin \bar{\Omega}} K^\alpha(z)dz, \quad \mbox{and} \quad 
\bar{\varphi}(x) = \int \limits_{x + z \notin \bar{\Omega}} \varphi(x + z) K^\alpha(z)dz,
\end{equation*}

With this, we remark that problem~\eqref{eq}-\eqref{Dirichlet} can be seen as an elliptic problem of the type
\begin{equation}\label{eqcensoredintro}
\bar{\lambda}(x) u(x) - \mathcal{I}_\Omega[u](x) + H(x, Du(x)) = \bar{\varphi}(x), \quad x \in \Omega,
\end{equation}
with generalized boundary condition $u = \varphi$ on $\partial \Omega$.

Note that the censored operator $\mathcal{I}_\Omega$ ``localizes'' inside $\Omega$ the equation, easing the analysis of the problem 
because, as we already explained in~\eqref{eqintro}, an equation like~\eqref{eqcensoredintro} concerns unknown functions defined in 
$\bar{\Omega}$, incorporating the exterior condition into the equation through the function $\bar{\varphi}$. However, by~\eqref{ellipticity} 
this function is potentially unbounded at the boundary and this creates a difficulty because we have to deal with unbounded data. 

The idea is to use a censored equation like~\eqref{eqcensoredintro}, taking advantage of the operator $\mathcal{I}_\Omega$ but
taking care about the unbounded term $\bar{\varphi}$.

\medskip

\noindent
{\bf Organization of the Paper.} In section~\ref{results} we introduce the principal assumptions and provide the main results. In 
section~\ref{notationsection} we state the basic notation, give the precise notion of 
solution of~\eqref{eq}-\eqref{Dirichlet} and introduce the censored problem. In section~\ref{technicalsection} we provide several technical 
results, including the key relation among our original problem and the censored one. These results allow us to undestand
the behavior of the sub and supersolutions at the boundary in section~\ref{propertiessection}, and an important improvement on its semicontinuity in 
section~\ref{coneconditionsection}. This last property will be the key fact to prove Theorems~\ref{teo1} and~\ref{teo2} in section~\ref{proofteo}. Finally, 
some extensions are discussed in section~\ref{remarkssection}.


\section{Assumptions and Main Results.}
\label{results}

In this section we complete the set of hypotheses allowing us to accomplish the comparison results. 
In addition to~\eqref{ellipticity}, we assume the following \textsl{nondegeneracy condition} relating the nonlocal term with
$\lambda$
\begin{equation}\tag{{\bf M}} \label{measure}
\exists \ \mu_0 > 0 \ : \ \lambda + \inf \limits_{x \in \Omega} \displaystyle{\int \limits_{x + z \notin \Omega}} K(z)dz \geq \mu_0.
\end{equation}

Recalling the definition of $H$ in~\eqref{Hamiltonian}, we also assume the uniform Lipschitz continuity of $b$
\begin{equation}\tag{{\bf L}}\label{L}
\begin{split}
(\exists \ L > 0) \ (\forall \ \beta \in \mathcal{B}) \ (\forall \ x,y \in \bar{\Omega}) \ : \  |b_\beta(x) - b_\beta(y)| \leq L|x - y|.
\end{split}
\end{equation}

Condition~\eqref{L} is rather classical and implies the existence and uniqueness of the SDE~\eqref{trajectories}. 
Concerning~\eqref{measure}, we note that it implies the associated pay-off for these jump processes (see~\eqref{payoff}) is well defined. 
In fact, if there is no discount rate ($\lambda = 0$), the resultant condition~\eqref{measure}
says that each point of $\Omega$ can jump outside with non-zero probability, concluding that the exit time $\tau$ relative to~\eqref{trajectories}
is  finite a.e.. 

Because of the weak diffusion setting of our problem, we take care about the behavior of the drift terms at the boundary dividing 
it into representative sets. First, we consider
\begin{eqnarray*}
\Gamma_{in} = \{ x \in \partial \Omega \ : \ \forall \beta \in \mathcal{B}, \ b_\beta(x) \cdot Dd(x) > 0 \},
\end{eqnarray*}
which can be understood as the sets of points where all the drift terms push inside $\Omega$ the trajectories defined in~\eqref{trajectories}.
For this, it is reasonable to expect a loss the boundary condition there. However, as we will 
see later on, it is possible to get an important improvement on the semicontinuity of sub and supersolutions on $\Gamma_{in}$, a key fact to 
get the comparison up to the boundary. 

We also consider the set
\begin{equation*}
\Gamma_{out} = \{ x \in \partial \Omega \ : \ \forall \beta \in \mathcal{B}, \ b_\beta(x) \cdot Dd(x) \leq 0 \},
\end{equation*}
where it is reasonable to think that there is no loss of the boundary condition. 

We also have in mind an intermediate situation, considering the set
\begin{eqnarray*}
\Gamma := \partial \Omega \setminus (\Gamma_{in} \cup \Gamma_{out}).
\end{eqnarray*}

In this set, loss of the boundary condition cannot be discarded. However, by the very definition of $\Gamma$, for each $x \in \Gamma$ there exists 
$\bar{\beta},\underline{\beta}\in \mathcal{B}$ such that
\begin{equation}\label{Gammaexplanation}
b_{\underline{\beta}}(x) \cdot Dd(x) \leq 0 < b_{\bar{\beta}}(x) \cdot Dd(x).
\end{equation}

Dropping the supremum in the definition of $H$ by taking $\bar{\beta}$ or $\underline{\beta}$, a subsolution
of~\eqref{eq}-\eqref{Dirichlet} is also a subsolution for the linear equation associated to these controls and therefore
it enjoys the good properties of both $\Gamma_{out}$ and $\Gamma_{in}$ at points on $\Gamma$.

Concerning these subsets of the boundary, we assume
\begin{equation}\tag{{\bf H}} \label{H}
\Gamma_{in}, \Gamma_{out} \ \mbox{and} \ \Gamma \ \mbox{are connected components of} \ \partial \Omega.
\end{equation}

By the smoothness of $\partial \Omega$, 
each of these subsets is uniformly away the others, avoiding two completely different drift's behavior for arbitrarily close points.

Since the uncertainty of the value of the sub and supersolutions on the boundary is one of the main difficulties in the 
study of problem~\eqref{eq}-\eqref{Dirichlet}, we introduce the following definition: For $u$ upper semicontinuous in $\bar{\Omega}$,
$v$ lower semicontinuous in $\bar{\Omega}$ (which will be thought as sub and supersolution, respectively) denote
\begin{equation}\label{deftildeuv}
\begin{split}
\tilde{u}(x) &= \left \{ \begin{array}{lccl} &\limsup \limits_{y \in \Omega, \ y \to x} 
u(y)& \quad & \mbox{if} \ x \in \Gamma \cup \Gamma_{in}\\
&u(x)& \quad & \mbox{if} \ x \in \Omega \cup \Gamma_{out}. \end{array} \right . \\
\tilde{v}(x) &= \left \{ \begin{array}{lccl} &\liminf \limits_{y \in \Omega, \ y \to x} 
v(y)& \quad & \mbox{if} \ x \in \Gamma_{in} \\
&v(x)& \quad & \mbox{if} \ x \in \Omega \cup \Gamma \cup \Gamma_{out}. \end{array} \right .
\end{split}
\end{equation}

Clearly $\tilde{u} \leq u$, $v \leq \tilde{v}$, $\tilde{u}$ is upper semicontinuous in $\bar{\Omega}$ and $\tilde{v}$ is lower semicontinuous 
in $\bar{\Omega}$. 

From now, we write usc for upper semicontinuous and lsc for lower semicontinuous functions. The main result of the article is the following

\begin{teo}\label{teo1}
Assume~\eqref{ellipticity},~\eqref{measure},~\eqref{H},~\eqref{L} hold. Let $u$ be a bounded usc viscosity subsolution of~\eqref{eq}-\eqref{Dirichlet}
and $v$ a bounded lsc viscosity supersolution of~\eqref{eq}-\eqref{Dirichlet}. Then
$$
u \leq v \quad \mbox{in} \ \Omega.
$$

Moreover, if we define $\tilde{u}, \tilde{v}$ as in~\eqref{deftildeuv}, then $\tilde{u} \leq \tilde{v}$ in $\bar{\Omega}$.
\end{teo}

By the presence of loss of the boundary condition, this comparison principle is established by the use of the notion of generalized
boundary conditions for viscosity sub and supersolution given by H. Ishii in~\cite{Ishii1} (see also \cite{usersguide}, $\S 7$). Once 
the comparison holds, the configuration of problem~\eqref{eq}-\eqref{Dirichlet} makes posible the use of Perron's method for 
integro-differential equations
(see~\cite{Alvarez-Tourin},~\cite{Barles-Imbert},~\cite{Sayah2} and~\cite{usersguide},~\cite{Ishii2} for an introduction on the method)
to get as a corollary the following

\begin{teo}\label{teo2}
Assume~\eqref{ellipticity},~\eqref{measure},~\eqref{H}, ~\eqref{L} hold. Then, there exists a unique viscosity solution 
to problem~\eqref{eq}-\eqref{Dirichlet} in $C(\bar{\Omega})$. 
\end{teo}


\section{Notation and Notion of Solution.}
\label{notationsection}

\subsection{Basic Notation.} For $\delta > 0$ and $x \in \R^n$ we write $B_\delta(x)$ as the open ball of radius $\delta$ centered at $x$ and 
$B_\delta$ if $x = 0$. For an arbitrary set $A$, we denote $d_A(x) = dist(x, \partial A)$ the signed distance function to $\partial A$ 
which is nonnegative for $x \in A$ and nonpositive for $x \notin A$. For $\Omega$ we simply
write $d(x) = d_{\partial \Omega}(x)$ and define the set $\Omega_\delta$ as the open set of all 
$x \in \Omega$ such that $d(x) < \delta$. By the smoothness assumption over the domain, there exists a fixed number $\delta_0 > 0$,
depending only on $\Omega$, such that $d$ is smooth in the set of points $x$ such that $|d(x)| < \delta_0$. We write for 
$x \in \R^n$ and $\lambda \in \R$ the sets
$$
\Omega - x = \{ z  : x + z \in \Omega\} \quad \mbox{and} \quad \lambda \Omega = \{ \lambda z : z \in \Omega\}. 
$$

Finally, for a function $\phi: \R^n \to \R$ and $x,z \in \R^n$, we denote
\begin{equation*}
\delta(\phi, x, z) = \phi(x + z) - \phi(x).
\end{equation*}

\subsection{Notion of Solution.}
As we mentioned in the introduction, we will
understand the solutions as functions defined in $\bar{\Omega}$, fixing their value as $\varphi$ outside $\Omega$. However, we need to define 
properly the value on the boundary in order to fit in the definition of viscosity solution in the literature. 

For an usc (resp. lsc) function $u : \bar{\Omega} \to \R$, we define its upper (resp. lower) $\varphi-$extension as 
\begin{equation}\label{Uextension}
u^\varphi(x) \ (\mbox{resp.} \  u_\varphi(x)) = 
\left \{ \begin{array}{ll} u(x) & \mbox{if} \ x \in \Omega \\ \varphi(x) & \mbox{if} \ x \in \bar{\Omega}^c \\ 
\max \ (\mbox{resp.} \ \min) \{ u(x), \varphi(x) \} & \mbox{if} \ x \in \partial \Omega,\end{array} \right . 
\end{equation}

Note that if $u$ is usc, then $u^\varphi$ is the usc envelope in $\R^n$ of the function 
$u \mathbf{1}_{\bar{\Omega}} + \varphi \mathbf{1}_{\bar{\Omega}^c}$. Analogously, if $u$ is lsc, $u_\varphi$ is the lsc envelope in $\R^n$ 
of the same function. 

In what follows we consider $x \in \bar{\Omega}$, $\delta > 0$ and for any $\phi \in C^1(\bar{B}_\delta(x))$ and 
for any bounded semicontinuous function $w$, we define the operators
\begin{equation*}
\begin{split}
\mathcal{I}_{\delta} [\phi](x) &= \int \limits_{B_\delta} \delta(\phi, x,z) K^\alpha(z)dz \\
\mathcal{I}^\delta [w](x) &= \int \limits_{B_\delta^c \cap (\Omega - x)} \delta(w, x,z) K^\alpha(z) dz
+ \int \limits_{B_\delta^c \cap (\Omega^c - x)} [\varphi(x + z) - w(x)] K^\alpha(z)dz.
\end{split}
\end{equation*}
and, for each $\beta \in \mathcal{B}$ 
\begin{equation}\label{Eevaluation}
E_{\beta, \delta}(w, \phi, x) = \lambda w(x) - \mathcal{I}_{\delta}[\phi](x) - \mathcal{I}^\delta[w](x) 
- b_\beta(x) \cdot D\phi(x) - f_\beta(x),
\end{equation}
where ``$E$'' stands for ``evaluation''. The nonlinear evaluation reads as
\begin{equation}\label{Enonlinear}
 E_{\delta}(w, \phi, x) = \sup \limits_{\beta \in \mathcal{B}} E_{\beta, \delta}(w, \phi, x).  
\end{equation}

The precise notion of solution used here is given by the following

\begin{defi}\label{defsolproblem}{\bf (Solution to Problem~\eqref{eq}-\eqref{Dirichlet})}

An usc function $u: \bar{\Omega} \to \R$ is a viscosity subsolution to~\eqref{eq}-\eqref{Dirichlet} if for each 
$x_0 \in \bar{\Omega}$ and each smooth function $\phi: \R^n \to \R$ such that $x_0$ is a maximum point of $u^\varphi - \phi$ in $B_\delta(x_0)$, then 
\begin{equation*}
\begin{split}
E_\delta(u^\varphi, \phi, x_0) \leq 0 & \quad \mbox{if} \  x_0 \in \Omega, \\
\min \{ E_\delta(u^\varphi, \phi, x_0), u^\varphi(x_0) - \varphi(x_0) \} \leq 0 & \quad \mbox{if} \ x_0 \in \partial \Omega.
\end{split}
\end{equation*}

A lsc function $v: \bar{\Omega} \to \R$ is a viscosity supersolution to~\eqref{eq}-~\eqref{Dirichlet} if for each 
$x_0 \in \bar{\Omega}$ and each smooth function $\phi: \R^n \to \R$ such that $x_0$ is a minimum point of $v_\varphi - \phi$ in $B_\delta(x_0)$, then
\begin{equation*}
\begin{split}
E_\delta(v_\varphi, \phi, x_0) \geq 0 & \quad \mbox{if} \  x_0 \in \Omega, \\
\min \{ E_\delta(v_\varphi, \phi, x_0), v_\varphi(x_0) - \varphi(x_0) \} \geq 0 & \quad \mbox{if} \ x_0 \in \partial \Omega.
\end{split}
\end{equation*}

Finally, a solution to~\eqref{eq}-\eqref{Dirichlet} is a function in $C(\bar{\Omega})$ which is simultaneously sub and supersolution in the above sense.
\end{defi}

The above definition is basically the one presented in~\cite{Barles-Imbert},~\cite{Barles-Chasseigne-Imbert},~\cite{Sayah1} 
and~\cite{Sayah2}. Written in that way we highlight the goal of this paper, which is to state the existence and uniqueness of a solution 
of~\eqref{eq}-\eqref{Dirichlet} in $C(\bar{\Omega})$. However, this definition is not the most confortable to deal with the discontinuities 
at the boundary of $u_\varphi$ and $u^\varphi$. For this reason in some situations we will look for the equation without
taking into account the exterior condition through the operator~\eqref{censoredI}, which ``censores'' the jumps outside $\Omega$. Hence, 
we need to introduce some notation in order to write properly the viscosity inequalities concerning this operator.

If $\phi$ is $C^1$ and bounded, and if $w$ is a bounded and semicontinuous function, the viscosity evaluation for this operator are defined by
\begin{equation}\label{nonlocalOmega}
\begin{split}
\mathcal{I}_{\Omega, \delta}[\phi](x) &= \int_{B_\delta \cap (\Omega - x)} \delta(\phi, x,z)K^\alpha(z)dz, \\
\mathcal{I}_{\Omega}^\delta[w](x) &= \int_{B_\delta^c \cap (\Omega - x)} \delta(w, x,z)K^\alpha(z)dz.
\end{split}
\end{equation}

The idea is to use this censored operator in auxiliary problems and therefore
distinct to equation~\eqref{eq}. To do so, we provide a viscosity notion of solution in a rather general setting. We consider a set
$A \subseteq \bar{\Omega}$, a relatively open subset of $\bar{\Omega}$ and $F: A \times \R \times \R^n \to \R$ a continuous function.

\begin{defi}\label{defcensored}{\bf (Solution for the Censored Equation)}

An usc function $u: \bar{\Omega} \to \R$ is a viscosity subsolution to
\begin{equation}\label{eqcensored}
-\mathcal{I}_\Omega[w] + F(x, u, Du) = 0 \quad \mbox{in} \ A
\end{equation}
if, for each $x_0 \in A$ and any smooth function $\phi: \R^n \to \R$ such that $x_0$ is a maximum point of $u - \phi$ in 
$A \cap B_\delta(x_0)$, then 
\begin{equation*}
-\mathcal{I}_\Omega^\delta[u](x_0) - \mathcal{I}_{\Omega, \delta}[\phi](x_0) + F(x_0, u(x_0), D\phi(x_0)) \leq 0.
\end{equation*}

A lsc function $v: \bar{\Omega} \to \R$ is a viscosity subsolution to~\eqref{eqcensored} if for each $x_0 \in A$ and each smooth function 
$\phi: \R^n \to \R$ such that $x_0$ is a minimum point of $v - \phi$ in $A \cap B_\delta(x_0)$, then 
\begin{equation*}
-\mathcal{I}_\Omega^\delta[v](x_0) - \mathcal{I}_{\Omega, \delta}[\phi](x_0) + F(x_0, v(x_0), D\phi(x_0)) \geq 0. 
\end{equation*}

A viscosity solution to~\eqref{eqcensored} is a function $u: \bar{\Omega} \to \R$, continuous in $\bar{\Omega}$ which is simultaneously sub and 
supersolution. 
\end{defi}

In what follows, when we say ``solution to~\eqref{eq}-\eqref{Dirichlet}'' we mean that it is in the sense of 
Definition~\ref{defsolproblem}. Otherwise, we mean it is a solution just to an equation and then we are refering to Definition~\ref{defcensored}.
The same applies for sub and supersolutions.

The following is technical notation to be used in this paper. We denote for $x,p \in \R^n$ 
\begin{equation}\label{HiHs}
\mathcal{H}_i(x, p) = \inf \limits_{\beta \in \mathcal{B}} \{ -b_\beta(x) \cdot p\}; \quad 
\mathcal{H}_s(x, p) = \sup \limits_{\beta \in \mathcal{B}} \{ -b_\beta(x) \cdot p\}.
\end{equation}

In order to depict the viscosity evaluation for equations associated to the above Hamiltonians and the censored operators
$\mathcal{I}_\Omega$ defined in~\eqref{nonlocalOmega}, we write for $x \in \bar{\Omega}$, $\beta \in \mathcal{B}$, $\phi$ smooth and $w$ bounded 
\begin{equation}\label{OmegaE}
\begin{split}
& \mathcal{E}_{\delta, \beta}(w, \phi, x) = - \mathcal{I}_{\Omega, \delta}[\phi](x) - \mathcal{I}_\Omega^\delta[w](x) 
- b_\beta(x) \cdot D\phi(x), \\
& \mathcal{E}_{\delta}^i(w, \phi, x) = \inf \limits_{\beta \in \mathcal{B}} \{ \mathcal{E}_{\delta, \beta}(x, w, \phi) \}
\quad \mbox{and} \quad \mathcal{E}_{\delta}^s( w, \phi, x) = \sup \limits_{\beta \in \mathcal{B}} \{ \mathcal{E}_{\delta, \beta}(x, w, \phi) \}.
\end{split}
\end{equation}


\section{Preliminary Technical Results.}
\label{technicalsection}

We start with the following result, which asserts that subsolutions of the problem~\eqref{eq}-\eqref{Dirichlet} in the sense of 
Definition~\ref{defsolproblem} are also subsolutions for an associated censored problem. This fact will be a useful tool to see 
problem~\eqref{eq}-\eqref{Dirichlet} in a subtle different perspective, ``localizing'' it in $\Omega$. 

\begin{lema}\label{lemaeqequivlocal}
Let $u$ be a bounded viscosity subsolution of problem~\eqref{eq}-\eqref{Dirichlet}. Then, there exists a constant $\theta_0 > 0$ 
such that $u$ is a viscosity subsolution to 
\begin{eqnarray}\label{eqequivlocal}
\lambda u - \mathcal{I}_\Omega[u] + H(x,Du) = \theta_0 (||\varphi||_\infty + ||u||_\infty) d^{-\alpha} \quad \mbox{in} \ \Omega.
\end{eqnarray}

Analogously, if $v$ is a bounded viscosity supersolution
of~\eqref{eq}-\eqref{Dirichlet}, then it is a viscosity supersolution to~\eqref{eqequivlocal} replacing $\theta_0$  by $-\theta_0$. 
The constant $\theta_0$ depends only on $n, \alpha$ and $\Lambda$. 
\end{lema}

\noindent
{\bf Proof:} Let $x \in \Omega$ and $\phi$ a smooth function such that $x$ is a maximum point of $u - \phi$ in $\Omega \cap B_\delta(x)$. 
For $\delta' \leq \min \{ \delta, d(x) \}$ we use that $u$ is a subsolution to~\eqref{eq} to write
\begin{equation*}
\lambda u(x) - \mathcal{I}_{\delta'}[\phi](x) - \mathcal{I}^{\delta'}[u^\varphi](x) + H(x, D\phi(x)) \leq 0.
\end{equation*}

Note that $\delta(u, x, z) \leq \delta(\phi, x, z)$ in $B_\delta \cap \Omega\setminus B_{\delta'}$. Applying this on the last inequality we arrive to
\begin{equation*}
\lambda u(x) - \mathcal{I}_{\Omega, \delta}[\phi](x) - \mathcal{I}_\Omega^\delta[u](x) + H(x, D\phi(x))
\leq (||\varphi||_\infty + ||u||_\infty) \int \limits_{B_{d(x)}^c} K^\alpha(z)dz,
\end{equation*}
where we have used notations~\eqref{nonlocalOmega}. Using~\eqref{ellipticity} we have
\begin{equation*}
\int \limits_{B_{d(x)}^c} K^\alpha(z)dz \leq \Lambda |\partial B_1| d(x)^{-\alpha} \int \limits_{1}^{+\infty}r^{-(1 + \alpha)}dr 
= \Lambda |\partial B_1| \alpha^{-1} d(x)^{-\alpha}, 
\end{equation*}
where $|\partial B_1|$ is the Lebesgue measure of the unit sphere in $\R^n$. The result follows taking $\theta_0 = \alpha^{-1} \Lambda |\partial B_1|$.
\qed

The following technical results introduce particular functions to handle the difficulties 
arising in the analysis on $\Gamma$ and $\Gamma_{in}$. As we will see, the assumption $\alpha < 1$ plays a crucial role to get them. 
The next lemma is intended to deal with the unpleasant term in the right hand side of~\eqref{eqequivlocal}. 
\begin{lema}\label{lemad}
Let $x_0 \in \partial \Omega$ and $b$ a Lipschitz continuous vector field. Assume there exists $c_0 > 0$ such that $b$ satisfies
\begin{eqnarray}\label{bpointinside0}
b(x_0) \cdot Dd(x_0) \geq c_0.
\end{eqnarray}

Then, for all $0 < \sigma$, there exists $0 < \bar{r} < \delta_0$\footnote{Recall that $\delta_0 > 0$ is such that $d$ 
is smooth in $\{ x \in \R^n : |d(x)| < \delta_0 \}$.} and $0 < \tilde{c}_0 < c_0$ such that the function $d^\sigma$ is a classical 
subsolution to the equation
\begin{eqnarray*}\label{eqd}
-\mathcal{I}_{\Omega}[w] - b \cdot Dw = -\tilde{c}_0 \ d^{\sigma - 1} \quad \mbox{in} \ B_{\bar{r}}(x_0)\cap \Omega.
\end{eqnarray*}
\end{lema}

\noindent
{\bf Proof:} First, we take care about the nonlocal term. We denote
\begin{equation}\label{Theta}
\Theta_x = ( \Omega - x ) \cap \{ z : d(x + z) < d(x) \}, 
\end{equation}
where $d = d_{\partial \Omega}$. In particular we have $\delta(d^\sigma, x, z) \geq 0$ for $z \in (\Omega - x) \setminus \Theta_x$ by the monotony 
of the function $t \mapsto t^\sigma$. Using this and condition~\eqref{ellipticity}, we get
\begin{eqnarray*}
-\mathcal{I}_{\Omega}[d^\sigma](x) \leq \Lambda \int \limits_{\Theta_x} -\delta(d^\sigma, x,z)|z|^{-(n + \alpha)}dz. 
\end{eqnarray*} 

The idea is to estimate the last integral splitting the domain $\Theta_x$. Simply dropping the nonpositive term $-d^\sigma(x + z)$, 
there exists $C > 0$ independent of $x$ such that
\begin{eqnarray*}
\begin{split}
-\int \limits_{\Theta_x \setminus B_{d(x)/2}} \delta(d^\sigma, x,z)|z|^{-(n + \alpha)}dz \leq & \ 
d(x)^\sigma \int \limits_{B_{d(x)/2}^c} |z|^{-(n + \alpha)}dz \\
= & \ 2^\alpha \alpha^{-1}|\partial B_1| d(x)^{\sigma - \alpha}
\end{split}
\end{eqnarray*}
meanwhile, using the differentiability of the distance function we have
$$
|\delta(d^\sigma, x, z)| \leq 2^{1 - \sigma} \sigma d^{\sigma - 1}(x)|z|,
$$
for all $z \in B_{d(x)/2}$. Thus, we conclude
\begin{eqnarray*}
\begin{split}
-\int \limits_{B_{d(x)/2}} \delta(d^\sigma, x, z)|z|^{-(n + \alpha)}dz \leq & \ 2^{1 - \sigma}|\partial B_1| \sigma d(x)^{\sigma - 1} 
\int \limits_{0}^{d(x)/2} r^{-\alpha} dr \\
= & \ 2^{\alpha - \sigma} (1 - \alpha)^{-1} \sigma |\partial B_1| d(x)^{\sigma - \alpha}.
\end{split}
\end{eqnarray*}

These last two estimates imply 
\begin{eqnarray}\label{estimatelemad}
-\mathcal{I}_{\Omega}[d^\sigma](x) \leq C d(x)^{\sigma - \alpha},
\end{eqnarray} 
for a constant depending only on the data. On the other hand, by~\eqref{bpointinside0} and the continuity of $x \mapsto b(x) \cdot Dd(x)$, 
we can assume 
$
b(x) \cdot Dd(x) \geq 3c_0/4,
$
for all $x \in B_{\bar{r}}(x_0) \cap \Omega$. Evaluating classically and using estimate~\eqref{estimatelemad}, we obtain
\begin{eqnarray*}
-\mathcal{I}_{\Omega}[d^\sigma](x) - b(x) \cdot Dd^\sigma(x) \leq d(x)^{\sigma-1} (-3 \sigma c_0/4  + Cd^{1 - \alpha}(x)).
\end{eqnarray*}

Since $\alpha < 1$, choosing $\bar{r}$ small depending on the data we get the result.
\qed

For each $r > 0$ small enough, define the set 
\begin{equation}\label{Sigma}
\Sigma_r = \{ x \in \Omega : 0 < dist(x, \Gamma_{in}) < r \}.   
\end{equation}

Using the compactness of $\mathcal{B}$ it is posible to get the following nonlinear version of the last lemma, valid in $\Gamma_{in}$.

\begin{cor}\label{cord}
There exists $\tilde{c}_0, \bar{r} > 0$ such that for all $\sigma > 0$, $d^\sigma$ is a classical subsolution to the equation
\begin{eqnarray*}\label{eqdnonlinear}
-\mathcal{I}_{\Omega}[w] + \mathcal{H}_s(x, Dw) = 
-\tilde{c}_0 \ d^{\sigma - 1} \quad \mbox{in} \ \Sigma_{\bar{r}},
\end{eqnarray*}
where $\mathcal{H}_s$ is defined in~\eqref{HiHs}.
\end{cor}

\noindent
{\bf Proof:} We remark that, by the compactness of $\mathcal{B}$, the continuity of $b$ and the smoothness of $\partial \Omega$, there exists a constant 
$c_0 > 0$ such that $b_\beta(x) \cdot Dd(x) \geq c_0$ for all $\beta \in \mathcal{B}$ and $x \in \Gamma_{in}$. Hence, we can obtain the result
of the previous proposition with the same constants $\bar{r}$ and $\tilde{c_0}$ for each point of $\Gamma_{in}$ and each control $\beta$. Taking
supremum over $\mathcal{B}$ we conclude the result.
\qed

For the next results we introduce the function $\zeta: \R^n \to \R$ defined as
\begin{eqnarray}\label{zeta}
\zeta(x) = \left \{ \begin{array}{cl} \log(d(x)) & \quad \mbox{if} \ x \in \Omega \\
 0 & \quad \mbox{if} \ x \in \Omega^c, \end{array} \right . 
\end{eqnarray}
which is bounded above, usc in $\R^n$, continuous in $\Omega$ and smooth in $\Omega_{\delta_0}$. Concerning this function, we need the following
estimate of its evaluation on the nonlocal operator. 

\begin{lema}\label{lemaIlog}
There exists constants $0 < \delta^* < \delta_0$ and $C > 0$ such that 
\begin{equation*}\label{eqIlog}
-\mathcal{I}_{\Omega}[\zeta]  \leq C d^{-\alpha} \quad \mbox{in} \ \Omega_{\delta^*}.
\end{equation*}

The constants $C, \delta^*$ depend only on the data and the smoothness of $\partial \Omega$.
\end{lema}

The proof of this estimate is very similar to the derivation of~\eqref{estimatelemad}, but we have to deal with the unboundedness of $\zeta$ so we 
postpone its proof to the Appendix.
Since $\zeta(x) \to -\infty$ as 
$x \to \partial \Omega$, the idea is to use $\zeta$ as a penalization for the testings of subsolutions in order to push the maximum test points
inside the domain. We can use properly this penalization argument proving $\zeta$ ``behaves'' like a subsolution for the problem when 
the drift term points strictly inside $\Omega$.

\begin{lema}\label{lematechnicallog}
Let $x_0 \in \partial \Omega$ and $b$ a Lipschitz continuous vector field. Assume there exists $c_0 > 0$ such that $b$ satisfies~\eqref{bpointinside0}.
Then, there exists $0 < \bar{r} < \delta_0$ and $0 < \tilde{c}_0 < c_0$, such that $\zeta$ is a classical subsolution to the equation
\begin{eqnarray}\label{eqtechnicallog}
-\mathcal{I}_{\Omega}[w] - b \cdot Dw = -\tilde{c}_0 \ d^{-1} \ \mbox{in} \ B_{\bar{r}}(x_0)\cap \Omega. 
\end{eqnarray}
\end{lema}

The proof of this result follows the same lines of Lemma~\ref{lemad} using Lemma~\ref{lemaIlog}. Again by the compactness of
$\mathcal{B}$, we have the following

\begin{cor}\label{corlog}
There exists $c_0, \bar{r} > 0$ such that $\zeta$ is a classical subsolution to the
equation
\begin{eqnarray*}\label{eqlognonlinear}
-\mathcal{I}_{\Omega}[w] + \mathcal{H}_s(x, Dw) = -c_0 \ d^{-1}  \quad \mbox{in} \ \Sigma_{\bar{r}},
\end{eqnarray*}
where $\mathcal{H}_s$ is defined in~\eqref{HiHs} and $\Sigma_r$ is defined in~\eqref{Sigma}.
\end{cor}

\begin{remark}\label{rmktechnical}
We can obtain the same results of Lemmas~\ref{lemad},~\ref{lemaIlog},~\ref{lematechnicallog} and its corollaries replacing $\mathcal{I}_\Omega$ 
by $\mathcal{I}_{\Omega'}$ for any $\Omega' \subseteq \Omega$. Moreover, the constant arising in each property does not change. 
\end{remark}


\section{Behavior of Sub and Supersolutions on the Boundary.}
\label{propertiessection}

\subsection{Classical Boundary Condition.}
Here we establish sufficient conditions to get the boundary condition in the classical sense. We start with the following result 
whose proof follows closely the arguments of~\cite{Barles-Chasseigne-Imbert}.

\begin{prop}\label{lemaGammaout}
Assume~\eqref{ellipticity},~\eqref{measure} and~\eqref{L} hold. Let $x_0 \in \partial \Omega$ and $u$, $v$ be bounded sub and supersolutions 
of problem~\eqref{eq}-\eqref{Dirichlet}, respectively. For $x \in \partial \Omega$, define 
\begin{equation*}
\mathcal{B}_{out}(x) = \{ \beta \in \mathcal{B} \ : \ b_\beta(x) \cdot Dd(x) \leq 0\}. 
\end{equation*}

\noindent
{\bf $(i)$} If there exists $r > 0$ such that $\mathcal{B}_{out}(x) \neq \emptyset$ for all $x \in B_{r}(x_0) \cap \partial \Omega$, then 
$$
u(x_0) \leq \varphi(x_0).
$$

\medskip
\noindent
{\bf $(ii)$} If there exists $r > 0$ such that $\mathcal{B}_{out}(x) = \mathcal{B}$ for all $x \in B_{r}(x_0) \cap \partial \Omega$, then
\begin{equation*}
\varphi(x_0) \leq v(x_0). 
\end{equation*}
\end{prop}

Note that, under the additional hypothesis~\eqref{H}, we have that each point $x \in \Gamma_{out}$ satisfies conditions $(i)$ and $(ii)$ in the 
above proposition, meanwhile each point $x \in \Gamma$ satisfies $(i)$. Thus, the inmediate consequence of Proposition~\ref{lemaGammaout} 
is the following

\begin{cor}\label{corlemaGammaout}
Assume~\eqref{ellipticity},~\eqref{measure},~\eqref{H} and~\eqref{L} hold. Let $u,v$ be bounded viscosity sub and supersolution 
to~\eqref{eq}-\eqref{Dirichlet}, respectively, and $\tilde{u}, \tilde{v}$ as in~\eqref{deftildeuv}. Then
\begin{equation*}
\begin{array}{cl} \tilde{u} \leq u \leq \varphi \leq v \leq \tilde{v}, \ & \mbox{in} \ \Gamma_{out},  \\ 
\tilde{u} \leq u \leq \varphi, \ & \mbox{in} \ \Gamma. \end{array}
\end{equation*}
\end{cor}

\noindent
{\bf Proof of Proposition~\ref{lemaGammaout}:} For $(i)$, we assume by contradiction that $u(x_0) - \varphi(x_0) = \nu$ for some $\nu > 0$, implying
that $u^\varphi(x_0) = u(x_0)$.
We consider $C^1$ functions $\chi, \psi: \R \to \R$ such 
that $\chi$ is even, bounded, $\chi(0) = 0$, $\chi(t) > 0$ for $t \neq 0$, $\liminf_{|t| \to \infty} \chi(t) > 0$ and such
that $\chi(t) = |t|^\sigma$ with $\sigma > 1$ in a neighborhood of $0$. For $\psi$ we assume it is
bounded, strictly increasing, $||\psi||_\infty \leq \frac{1}{4} \nu$ and $\psi(t) = k t$ in a neighborhood of $0$, for some $k > 0$. 

Consider a parameter $\eta > 0$ to be sent to zero and $\epsilon = \epsilon(\eta)$ such that $\epsilon(\eta) \to 0$ as $\eta \to 0$
in a rate to be precised later on. We introduce the following penalization function
$$
\Psi(y) := \epsilon^{-1} \chi(|y - x_0|) + \psi(d(y)/\eta).
$$

By definition of $u^\varphi$, the function $x \mapsto u^\varphi(x) - \Psi(x)$ is upper semicontinuous in $\R^n$. Using that
$\liminf_{|t| \to \infty} \chi(t) > 0$, taking $\epsilon$ suitably small we conclude $u^\varphi - \Psi$ attains its global maximum at a 
point $\bar{x} \in \R^n$. Using that 
\begin{equation}\label{ineqcorlemaGammaout}
u^\varphi(\bar{x}) - \Psi(\bar{x}) \geq u^\varphi(x_0) - \Psi(x_0) = u^\varphi(x_0), 
\end{equation}
we conclude the term $\Psi(\bar{x})$ is bounded for all $\eta$ and then $\bar{x} \to x_0$ as $\eta \to 0$. Now,
since $\nu + \varphi(x_0) = u^\varphi(x_0)$, we conclude 
$
\nu + \varphi(x_0) \leq u^\varphi(\bar{x}) - 1/4 \nu, 
$
meaning that $\bar{x} \in \bar{\Omega}$ for each $\eta$ small. Moreover, using again~\eqref{ineqcorlemaGammaout} and the upper semicontinuity 
of $u$, we have 
\begin{equation}\label{proplemaGammaout}
d(\bar{x}) = o_1(\eta) \eta, \ |\bar{x} - x_0| = o_\eta(1), \  \mbox{and} \ u^\varphi(\bar{x}) \to u(x_0),
\end{equation}
as $\eta \to 0$. Since we can use the penalization as a test function for $u$ at $\bar{x}$, for each $\beta \in \mathcal{B}$ and $\delta > 0$ we have
\begin{equation}\label{testingproplemaGammaout}
\lambda u^\varphi(\bar{x}) \leq \mathcal{I}_{\delta, \beta}[\Psi](\bar{x}) + \mathcal{I}_\beta^\delta[u](\bar{x}) + b_\beta(\bar{x}) \cdot D\Psi(\bar{x}) 
+ f_\beta(\bar{x}). 
\end{equation}

We need to estimate properly the nonlocal evaluations of the testing. For this, consider $r > 0$ small, but independent of $\eta$. Considering
$0 < \delta \leq \mu < r$ we define the sets
\begin{equation*}
\begin{split}
& \mathcal{A}_{\delta}^{ext} = \{ z \in B_r \ : \ d(\bar{x} + z) \leq d(\bar{x}) - \delta \}. \\ 
& \mathcal{A}_{\delta, \mu} = \{ z \in B_r \ : \ d(\bar{x}) - \delta < d(\bar{x} + z) < d(\bar{x}) + \mu \}. \\ 
& \mathcal{A}_{\mu}^{int} = \{ z \in B_r \ : \ \mu + d(\bar{x}) \leq d(\bar{x} + z) \}. 
\end{split}
\end{equation*}

We remark that $B_\delta \subset \mathcal{A}_{\delta, \mu}$ and that the constant $C > 0$ arising in each of the following estimates does not depend 
on $r, \mu, \delta, \eta$ or $\epsilon$. By~\eqref{ellipticity}, we clearly have 
\begin{equation*}
\begin{split}
& \int \limits_{B_r^c} \delta(u^\varphi, \bar{x}, z)K^\alpha(z)dz \leq C||u^\varphi||_\infty r^{-\alpha}, \\
& \int \limits_{\mathcal{A}_{\mu}^{int}} \delta(u^\varphi, \bar{x}, z)K^\alpha(z)dz \leq C||u||_\infty \mu^{-\alpha}.
\end{split}
\end{equation*}

By the continuity of $\varphi$ and the last fact in~\eqref{proplemaGammaout}, for all $r$ and $\eta$ small we have 
$\varphi(\bar{x} + z) \leq \varphi(x_0) + \nu/4$ for all $z \in \mathcal{A}_\delta^{ext}$ and $u(x_0) - \nu/4 \leq u^\varphi(\bar{x})$.
Then, since $u(x_0) = \varphi(x_0) + \nu$ and using~\eqref{ellipticity} we conclude
\begin{equation*}
\begin{split}
\int \limits_{\mathcal{A}_\delta^{ext}} \delta(u^\varphi, \bar{x}, z)K^\alpha(z)dz \leq
 -\nu/2 \int \limits_{\mathcal{A}_\delta^{ext}} K^\alpha(z)dz \leq -C \nu \delta^{-\alpha}, 
\end{split}
\end{equation*}
where $C$ depends on the data and $r$, but not on $\eta$. Finally, using that $\bar{x}$ is a global maximum point of $u - \Psi$, we
have $\delta(u^\varphi, \bar{x}, z) \leq \delta(\Psi, \bar{x}, z)$ and then since $\Psi$ is $C^1$ we get
\begin{equation}\label{estimatelemaGammaoutapp}
\begin{split}
\int \limits_{\mathcal{A}_{\delta, \mu}} \delta(u^\varphi,\bar{x}, z)K^\alpha(z)dz \leq &
\int \limits_{\mathcal{A}_{\delta, \mu}} \delta(\Psi,\bar{x}, z)K^\alpha(z)dz   \\
\leq & C(\eta^{-1} + \epsilon^{-1})\mu^{1 - \alpha}. 
\end{split}
\end{equation}

A detailed proof of the second inequality in~\eqref{estimatelemaGammaoutapp} can be found in the Appendix. Setting $\mu = \eta$, $\epsilon \geq \eta$
and $d(\bar{x}) < \delta \leq o_\eta(1)\eta$, we conclude from the last estimates that the nonlocal terms
in~\eqref{testingproplemaGammaout} can be computed as
\begin{equation*}
\mathcal{I}_{\delta, \beta}[\Psi](\bar{x}) + \mathcal{I}_\beta^\delta[u^\varphi](\bar{x}) \leq -(o_\eta(1)\eta)^{-\alpha}, 
\end{equation*}
so, replacing this in~\eqref{testingproplemaGammaout} and specifying the choice $\epsilon = \eta^\alpha$, we conclude
\begin{equation}\label{ineqLemaGammaout0}
\lambda u(\bar{x}) - f_\beta(\bar{x}) \leq k \eta^{-1} b_\beta(\bar{x}) \cdot Dd(\bar{x}) + C \eta^{-\alpha} |\bar{x} - x_0|^{\sigma - 1}
 - (o_\eta(1) \eta)^{-\alpha},
\end{equation}
for all $\beta \in \mathcal{B}$. Now, denoting $\hat{\bar{x}} \in \partial \Omega$ such that $d(\bar{x}) = |\bar{x} - \hat{\bar{x}}|$,
for all $\eta$ small there exists $\beta_\eta \in \mathcal{B}_{out}(\hat{\bar{x}})$, by the assumption in $(i)$. Thus, using~\eqref{L} and
the smoothness of the domain, there exists $C > 0$ independent of $\eta$ such that
\begin{eqnarray*}
b_{\beta_\eta}(\bar{x}) \cdot Dd(\bar{x}) \leq b_{\beta_\eta}(\bar{x}) \cdot Dd(\bar{x}) - b_{\beta_\eta}(\hat{\bar{x}}) \cdot Dd(\hat{\bar{x}})
\leq K|\bar{x} - \hat{\bar{x}}| = Cd(\bar{x}),
\end{eqnarray*}
and by~\eqref{proplemaGammaout} we conclude
\begin{equation*}
b_{\beta_\eta}(\bar{x}) \cdot Dd(\bar{x}) \leq o_\eta(1) \eta. 
\end{equation*}

Taking $\beta = \beta_\eta$ in~\eqref{ineqLemaGammaout0}, using the last inequality and~\eqref{proplemaGammaout} we conclude
\begin{equation}\label{ineqLemaGammaout}
\lambda u(\bar{x}) - ||f||_\infty \leq  o_\eta(1) + o_\eta(1) \eta^{-\alpha} - (o_\eta(1) \eta)^{-\alpha}.
\end{equation}

Letting $\eta \to 0$ we arrive to a contradiction because of the last term in the above inequality. 

For $(ii)$, we proceed in the same way using a similar
contradiction argument, reversing the signs of the penalizations in order to get a minimum test point for $v$. The role of $\beta_\eta$ above can 
be played by any $\beta \in \mathcal{B}$ by the assumption in $(ii)$ and then the estimates are independent of $\beta$ using ~\eqref{L}. 
Hence, we arrive to an expression very similar to~\eqref{ineqLemaGammaout}, with the reverse inequality, $v$ in place of $u$ 
and a plus sign in the last term in the right hand side. The contradiction is again obtained by letting $\eta \to 0$.
\qed

\begin{remark}
The lack of information concerning $v$ on $\Gamma$ has to do with the form of $H$: we simply cannot drop the supremum taking 
$\underline{\beta}$ given in~\eqref{Gammaexplanation} and ensure $v$ is still a supersolution. In fact, both cases $v \geq \varphi$ or $v < \varphi$
can happen depending on $f$ and $\varphi$.
\end{remark}


\subsection{Viscosity Inequality Up to the Boundary.}
Now we consider auxiliary functions associated to sub and supersolutions of problem~\eqref{eq}-\eqref{Dirichlet} which satisfy
the viscosity inequality up to the boundary for auxiliary censored problems. 

\begin{prop}\label{propU}
Let $x_0 \in \partial \Omega$ and $\beta_0 \in \mathcal{B}$ such that $b = b_{\beta_0}$ satisfies condition~\eqref{bpointinside0} for  
some $c_0 > 0$. Let $u$ be a bounded viscosity subsolution to~\eqref{eq}-\eqref{Dirichlet} and
$\tilde{u}$ as in~\eqref{deftildeuv}. Then, there exists $A, a > 0$ such that the function $U: \bar{\Omega} \to \R$ defined as
\begin{equation}\label{U}
U(x) = \tilde{u}(x) + A  d^{1 - \alpha}(x)
\end{equation}
is a viscosity subsolution to the equation
\begin{equation*}\label{eqpropU}	
- \mathcal{I}_\Omega[w] - b \cdot Dw = 0 \quad \mbox{in} \ B_a(x_0) \cap \bar{\Omega}.
\end{equation*}
\end{prop}

\noindent
{\bf Proof:} Since~\eqref{bpointinside0} holds, we can first choose $a < \bar{r}$ of Lemma~\ref{lemad}, which depends only on $c_0$, but not on $x_0$. 
Let $\delta > 0$, $\bar{x} \in B_a(x_0) \cap \Omega$ and $\phi$ a smooth function such that $\bar{x}$ is strict maximum point of 
$U - \phi$ in $B_\delta(x_0) \cap \bar{\Omega}$. Hence, we can see $\bar{x}$ as a test point for $u$ with test function $\phi - Ad^{1-\alpha}$ 
by the smoothness of $d$ near the boundary. Applying Lemma~\ref{lemaeqequivlocal} and the definition of $U$ we get
\begin{equation*}\label{ineqpropU}
\begin{split}
&- \mathcal{I}_{\Omega, \delta}[\phi](\bar{x}) - \mathcal{I}_{\Omega}^\delta[U](\bar{x}) - 
b(\bar{x}) \cdot D\phi(\bar{x}) \\
\leq & \ \lambda ||u||_\infty + f_{\beta_0}(\bar{x})
+ \theta_0 (||u||_\infty + ||\varphi||_\infty) d(\bar{x})^{-\alpha} \\
& + A(-\mathcal{I}_\Omega[d^{1 - \alpha}](\bar{x}) - b(\bar{x}) \cdot Dd^{1- \alpha}(\bar{x})), 
\end{split}
\end{equation*}
where we have used $\tilde{u} = u$ in $\Omega$. Using Lemma~\ref{lemad} in the last inequality we conclude that
there exists a constant $\tilde{c}_0 > 0$ such that
\begin{equation*}
\mathcal{E}_{\beta_0, \delta}(U, \phi, \bar{x}) \leq 
d^{-\alpha}(\bar{x}) \Big{(} Cd^\alpha(\bar{x}) + \theta_0 (||u||_\infty + ||\varphi||_\infty )- A\tilde{c}_0 \Big{)},
\end{equation*}
where we have used notation~\eqref{OmegaE}. Choosing $A > 0$ large enough and $a$ small depending on 
the data but not on $\bar{x}$ we arrive to
\begin{equation*}\label{ineqpropU2}
\mathcal{E}_{\beta_0, \delta}(U, \phi, \bar{x}) \leq -\tilde{c}_0/2 \ d(\bar{x})^{-\alpha}, 
\end{equation*}
concluding the proof in this case.

When $\bar{x} \in B_a(x_0) \cap \partial \Omega$, by definition of $\tilde{u}$ we consider a sequence $(x_k)_k$ of points in $\Omega$ such that 
$x_k \to \bar{x}$ and $u(x_k) \to \tilde{u}(\bar{x})$ and define $\epsilon_k = d(x_k)$. Let $\phi$ be a test function for 
$U$ in $B_\delta(\bar{x}) \cap \bar{\Omega}$ at $\bar{x}$  and consider the penalization
\begin{equation}\label{penalizationpropU}
x \mapsto \Phi(x) := U(x) - (\phi(x) - \epsilon_k \zeta(x)).
\end{equation}

From this, it is easy to see that  for all $k$ large there exists $\bar{x}_k \in \Omega$, maximum point of $\Phi$ in 
$B_\delta(\bar{x}) \cap \Omega$ and using the inequality $\Phi(\bar{x}_k) \geq \Phi(x_k)$ we get
\begin{equation}\label{proppropU}
\bar{x}_k \to \bar{x}, \ u(\bar{x}_k) \to \tilde{u}(\bar{x}), \quad \mbox{as} \ k \to \infty.
\end{equation}

Now we want to use the penalization as a test function, but since $\zeta$ is unbounded close to the boundary we have to restrict the testing set,
considering the last penalization as a testing in $B_{\delta'}(\bar{x}_k)$ with $\delta' < d(\bar{x}_k)$. Hence, since $\bar{x}_k \in \Omega$, 
arguing as above we conclude that
\begin{equation}\label{testingpropUforlog}
\mathcal{E}_{\beta_0, \delta'}(U, \phi, \bar{x}_k) 
\leq   -\tilde{c}_0/2 d(\bar{x}_k)^{-\alpha} + 
\epsilon_k \Big{(} -\mathcal{I}_{B_{\delta'}(\bar{x}_k)}[\zeta](\bar{x}_k) - b(\bar{x}_k) \cdot D\zeta(\bar{x}_k) \Big{)}.
\end{equation}

However, since we have $\bar{x}_k$ is a maximum point of the testing in $B_\delta(\bar{x}) \cap \Omega$, the inequality
\begin{equation*}
\delta(U, \bar{x}, z) \leq \delta(\phi, \bar{x}, z) - \epsilon_k \delta(\zeta, \bar{x}, z) \quad \mbox{for all} \ z \in (\Omega - \bar{x})
 \cap B_\delta
\end{equation*}
lead us to
\begin{equation*}
\begin{split}
& -\mathcal{I}_{\Omega}^\delta[U](\bar{x}_k) - \mathcal{I}_{\Omega, \delta}[\phi](\bar{x}_k) 
+ \epsilon_k \int \limits_{B_\delta \setminus B_{\delta'} \cap (\Omega - \bar{x}_k)} \delta(\zeta, \bar{x}_k, z) K^\alpha(z)dz \\
\leq & -\mathcal{I}_{\Omega}^{\delta'}[U](\bar{x}_k) - \mathcal{I}_{\Omega, \delta'}[\phi](\bar{x}_k).
\end{split}
\end{equation*}

Using this inequality in~\eqref{testingpropUforlog} we obtain the inequality
\begin{equation*}
\mathcal{E}_{\beta_0, \delta}(U, \phi, \bar{x}_k) 
\leq   -\tilde{c}_0/2 d(\bar{x}_k)^{-\alpha} + 
\epsilon_k \Big{(} -\mathcal{I}_{\Omega \cap B_\delta(\bar{x}_k)}[\zeta](\bar{x}_k) - b(\bar{x}_k) \cdot D\zeta(\bar{x}_k) \Big{)},
\end{equation*}
concluding, by Lemma~\ref{lematechnicallog} (see also Remark~\ref{rmktechnical}) that
\begin{equation*}
\mathcal{E}_{\beta_0, \delta}(U, \phi, \bar{x}_k) 
\leq   -\tilde{c}_0/2 d(\bar{x}_k)^{-\alpha} - \tilde{c}_0 \epsilon_k d(\bar{x}_k)^{-1} < 0.
\end{equation*}

Letting $k \to \infty$, using the smoothness of $\phi$, the continuity of $d$, and~\eqref{proppropU} together with the upper
semicontinuity  of $U$ in $\bar{\Omega}$ we get the result.
\qed

\begin{remark}
Note that the role of $\beta_0$ in the last proposition can be played by any control when $x_0 \in \Gamma_{in}$ and by $\bar{\beta}$ when
$x_0 \in \Gamma$, with $\bar{\beta}$ as in~\eqref{Gammaexplanation}.
\end{remark}

For the next result, we recall that $\Sigma_r = \{ x \in \Omega : dist(x, \Gamma_{in}) < r\}$.

\begin{prop}\label{propUV}
Let $u, v$ be respectively bounded viscosity sub and supersolution to problem~\eqref{eq}-\eqref{Dirichlet} and $\tilde{u}, \tilde{v}$ as 
in~\eqref{deftildeuv}. Recall $U$ given in~\eqref{U} and consider the function $V: \bar{\Omega} \to \R$ defined as 
\begin{equation}\label{V}
V(x) = \tilde{v}(x) - A  d^{1 - \alpha}(x).
\end{equation}

Then, there exist $A, a > 0$ such that the function $U$ (resp. $V$) is a viscosity subsolution (resp. supersolution) to the equation
\begin{equation*}
- \mathcal{I}_\Omega[w] + \mathcal{H}_s(x, Dw) = 0 \quad \mbox{in} \ \bar{\Sigma}_a.
\end{equation*}
\end{prop}

\noindent
{\bf Proof:} As in the proof of Corollary~\ref{cord}, we have the existence of a constant 
$c_0 > 0$ such that $b_\beta(x) \cdot Dd(x) \geq c_0$ for all $\beta \in \mathcal{B}$ and all $x \in \Gamma_{in}$.

Let $\bar{x} \in \Sigma_a$. If $\phi$ is a smooth function such that $\bar{x}$ is a maximum point for $U - \phi$ 
in $B_\delta(\bar{x}) \cap \bar{\Omega}$, 
proceeding as in Proposition~\ref{propU} and using notation~\eqref{OmegaE} we conclude that
\begin{equation*}
\begin{split}
\mathcal{E}_\delta^s(U, \phi, \bar{x}) 
\leq & \ \lambda||u||_\infty + \theta_0 (||u||_\varphi + ||\varphi||_\infty)d(\bar{x})^{-\alpha} + ||f||_\infty \\
& + A \Big{(} -\mathcal{I}_\Omega[d^{1 - \alpha}](\bar{x}) + \mathcal{H}_s(\bar{x}, Dd^{1 - \alpha}(\bar{x})) \Big{)},
\end{split}
\end{equation*}
and by the application of Corollary~\ref{cord}, taking $A$ large and $a$ small this inequality leads us to
\begin{equation}\label{testingpropUV}
\mathcal{E}_\delta^s(U, \phi, \bar{x}) \leq -\tilde{c}_0/2 \ d^{-\alpha}(\bar{x}), 
\end{equation}
concluding the proof of the proposition for this case. We highlight the constants $a$ and $A$ depend only 
on $c_0$ and not on the particular point $\bar{x}$ considered, concluding the result for $\Sigma_a$. To get the result in $\bar{\Sigma}_a$ we proceed
in the same way as in Proposition~\ref{propU}, considering $\bar{x} \in \partial \Omega$ and penalizing the testing by the introduction of the 
function $\epsilon_k \zeta$ in order to push the testing point inside $\Omega$ (see~\eqref{penalizationpropU} and its subsequent arguments).
Hence, using~\eqref{testingpropUV} we arrive to 
\begin{equation*}
\mathcal{E}_\delta^s(U, \phi, \bar{x}_k) \leq -\tilde{c}_0/2 \ d^{-\alpha}(\bar{x}_k)    
+ \epsilon_k \Big{(} -\mathcal{I}_{\Omega \cap B_\delta(\bar{x}_k)}[\zeta](\bar{x}_k) + \mathcal{H}_s(\bar{x}_k, D\zeta(\bar{x}_k)) \Big{)},
\end{equation*}
and applying Corollary~\ref{corlog} in the last expression we conclude 
\begin{equation}\label{testingpropUV2}
\mathcal{E}_\delta^s(U, \phi, \bar{x}_k) \leq -\tilde{c}_0/2 \ d^{-\alpha}(\bar{x}_k) - \tilde{c}_0 \epsilon_k d^{-1}(\bar{x}_k). 
\end{equation}

Thus, the right-hand side is nonpositive for all $k$ large and we conclude as 
in the proof of Proposition~\ref{propU} letting $k \to \infty$ by~\eqref{proppropU}. The result for $V$ can be obtained in the same way.
\qed

We finish this section with an analogous result for $U - V$. We shall introduce the following notation: For points 
$\bar{x}, \bar{y} \in \bar{\Omega}$ we denote
\begin{equation}\label{defD's}
\begin{split}
D_{int} = (\Omega - \bar{x}) \cap (\Omega - \bar{y}), \quad & D_{ext} =  (\Omega - \bar{x})^c \cap (\Omega - \bar{y})^c, \\ 
D_{int}^{\bar{x}} = (\Omega - \bar{x}) \cap (\Omega - \bar{y})^c, \quad & D_{int}^{\bar{y}} = (\Omega - \bar{x})^c \cap (\Omega - \bar{y}).
\end{split}
\end{equation}

We highlight that each set depends both on $\bar{x}$ and $\bar{y}$, but we omit the dependence on $\bar{x}$ and/or $\bar{y}$ in some cases for 
a sake of simplicity.

\begin{prop}\label{propW}
Let $u, v$  be respectively bounded viscosity sub and supersolution to problem~\eqref{eq}-\eqref{Dirichlet}, $U$ as in~\eqref{U}, $V$
as in~\eqref{V} and define the function $W = U - V$. Then, there exists $A, a > 0$ such that $W$ is a viscosity subsolution to the problem
\begin{equation*}
- \mathcal{I}_\Omega[w] + \mathcal{H}_i(x, Dw) = 0 \quad \mbox{in} \ \bar{\Sigma}_a,
\end{equation*}
where $\Sigma_r$ is defined in~\eqref{Sigma}.
\end{prop}

\noindent
{\bf Proof:} We start with the case the test point $x_0 \in \Omega$. Assume there exists $\phi$ a smooth function such that $x_0$ 
is a strict maximum point of $W - \phi$ in $B_\delta(x_0) \cap \bar{\Omega}$. For $\epsilon > 0$, consider the function 
$\phi_\epsilon(x,y) = \phi(x) + |x - y|^2/\epsilon^2$ and the penalization
\begin{equation}\label{penalizationcorpropUV}
(x,y) \mapsto \Phi(x,y) := U(x) - V(y) - \phi_\epsilon(x, y).
\end{equation}

This function is usc in $(B_\delta(x_0) \cap \bar{\Omega}) \times (B_\delta(x_0) \cap \bar{\Omega})$ and it attains its maximum at a 
point $(\bar{x}, \bar{y})$. Using the inequality $\Phi(\bar{x}, \bar{y}) \geq \Phi(x_0, x_0)$ we obtain as $\epsilon \to 0$
\begin{equation}\label{propcorUV}
|\bar{x} - \bar{y}|^2/\epsilon^2 \to 0, \quad 
\bar{x}, \bar{y} \to x_0, \quad \ U(\bar{x}) \to U(x_0), \quad V(\bar{y}) \to V(x_0).  
\end{equation}

Since $x_0 \in \Omega$, we have that for all $\epsilon$ suitably small,
$d(\bar{x}), d(\bar{y}) > d(x_0)/2$. Now we apply Proposition~\ref{propUV} looking the last penalization as a testing for $U$ and $V$. Using 
notation~\eqref{OmegaE}, we write for each $0 < \delta' < \delta$
\begin{equation*}
\mathcal{E}_{\delta'}^s(\bar{x}, U, \phi_\epsilon(\cdot, \bar{y})) \leq 0 \quad \mbox{and} \quad
\mathcal{E}_{\delta'}^s(\bar{y}, V, -\phi_\epsilon(\bar{x}, \cdot)) \geq 0
\end{equation*}
and then, for each $h > 0$ there exists $\beta_h \in \mathcal{B}$ such that
\begin{equation*}
\mathcal{E}_{\delta', \beta_h}(\bar{x}, U, \phi_\epsilon(\cdot, \bar{y})) \leq 0 \quad \mbox{and} \quad
\mathcal{E}_{\delta', \beta_h}(\bar{y}, V, -\phi_\epsilon(\bar{x}, \cdot)) \geq -h.
\end{equation*}

Substracting both inequalities, we are lead to
\begin{equation}\label{testingcorpropUV}
\begin{split}
& - \mathcal{I}_{\Omega, \delta'}[\phi](\bar{x}) - (\mathcal{I}_\Omega^{\delta'}[U](\bar{x}) - \mathcal{I}_\Omega^{\delta'}[V](\bar{y})) 
- b_{\beta_h}(\bar{x}) \cdot Dd(\bar{x}) \\
\leq & \ h + \epsilon^{-2} O({\delta'}^{1 - \alpha}) + 2 \epsilon^{-2}(b_{\beta_h}(\bar{x}) - b_{\beta_h}(\bar{x})) \cdot (\bar{x} - \bar{y}).
\end{split}
\end{equation}

By ~\eqref{L} and the first fact in~\eqref{propcorUV} we have the estimate
\begin{equation}\label{driftcorpropUV}
|(b_{\beta_h}(\bar{x}) - b_{\beta_h}(\bar{y}))\cdot (\bar{x} - \bar{y})/\epsilon^2| \leq o_\epsilon(1).
\end{equation}

In what follows we deal with the nonlocal terms in the left-hand side of~\eqref{testingcorpropUV}. 
Using notation~\eqref{defD's} we have 
\begin{equation*}
\Omega - \bar{x} = D_{int} \cup D_{int}^{\bar{x}}, \quad \mbox{and} \quad \Omega - \bar{y} = D_{int} \cup D_{int}^{\bar{y}}.
\end{equation*}

Thus, we write
\begin{equation}\label{nonlocalcorpropUV0}
\mathcal{I}_\Omega^{\delta'}[U](\bar{x}) - \mathcal{I}_\Omega^{\delta'}[V](\bar{y}) = J_{int, \bar{x}}^{\delta'} + J_{int, \bar{x}}^{\delta'} 
+ J_{int}^{\delta'},
\end{equation}
where 
\begin{equation*}
\begin{split}
& J_{int, \bar{x}}^{\delta'} = \int \limits_{D_{int}^{\bar{x}} \setminus B_{\delta'}}\delta(U,\bar{x}, z)K^\alpha(z)dz; \quad
J_{int, \bar{y}}^{\delta'} = \int \limits_{D_{int}^{\bar{y}} \setminus B_{\delta'}}\delta(V,\bar{x}, z)K^\alpha(z)dz; \\
& J_{int}^{\delta'} = \int \limits_{D_{int} \setminus B_{\delta'}}[U(\bar{x} + z) - V(\bar{y} + z) - (U(\bar{x}) - V(\bar{y}))]K^\alpha(z)dz.
\end{split}
\end{equation*}

By~\eqref{propcorUV}, as $\epsilon \to 0$ we have $\bar{x}, \bar{y} \to x_0$ and $|D_{int}^{\bar{x}}|, |D_{int}^{\bar{y}}| \to 0$. Then
\begin{equation}\label{nonlocalcorpropUV}
J_{int, \bar{x}}^{\delta'}, \ J_{int, \bar{y}}^{\delta'} = o_\epsilon(1) d(x_0)^{-\alpha}.
\end{equation}
with $o_\epsilon(1)$ independent of $\delta$. Meanwhile, since $(\bar{x}, \bar{y})$ is a maximum point, for all 
$z \in D_{int} \cap B_\delta \setminus B_{\delta'}$ we have
\begin{equation*}
U(\bar{x} + z) - V(\bar{y} + z) - (U(\bar{x}) - V(\bar{y})) \leq \delta(\phi, \bar{x}, z),
\end{equation*}
and then we obtain that
\begin{equation*}
\mathcal{I}_{\Omega, \delta'}[\phi](\bar{x}) + J_{int}^{\delta'} \leq \mathcal{I}_{\Omega, \delta}[\phi](\bar{x}) + J_{int}^{\delta}
+ \int \limits_{D_{int}^{\bar{x}} \cap B_\delta \setminus B_{\delta'}} \delta(\phi, \bar{x}, z)K^\alpha(z)dz,
\end{equation*}
but the last term in the right-hand side of the last inequality is $o_\epsilon(1)$ by the smoothness of $\phi$ and the fact $D_{int}^{\bar{x}}$
vanishes as $\epsilon \to 0$. Thus, we conclude that
\begin{equation}\label{delta'bydelta}
\mathcal{I}_{\Omega, \delta'}[\phi](\bar{x}) + J_{int}^{\delta'} \leq \mathcal{I}_{\Omega, \delta}[\phi](\bar{x}) + J_{int}^{\delta} + o_\epsilon(1). 
\end{equation}

Replacing this last expression together with~\eqref{nonlocalcorpropUV} in~\eqref{nonlocalcorpropUV0} we conclude
\begin{equation*}
\mathcal{I}_{\Omega, \delta'}[\phi](\bar{x}) + \mathcal{I}_\Omega^{\delta'}[U](\bar{x}) - \mathcal{I}_\Omega^{\delta'}[V](\bar{y}) 
\leq  o_\epsilon(1) d(x_0)^{-\alpha} + \mathcal{I}_{\Omega, \delta}[\phi(\bar{x})] + J_{int}^{\delta}.
\end{equation*}

Finally, replacing this last inequality together with~\eqref{driftcorpropUV} in~\eqref{testingcorpropUV} and taking infimum over 
$\beta_h \in \mathcal{B}$ in the right-hand side of~\eqref{testingcorpropUV}, we arrive to
\begin{equation*}\label{testingcorpropUVI}
\begin{split}
&- \mathcal{I}_{\Omega, \delta} [\phi](\bar{x}) - J_{int}^\delta + \mathcal{H}_i(\bar{x}, D\phi(\bar{x})) \\
\leq & \ h + \epsilon^{-2} O({\delta'}^{1 - \alpha}) + o_\epsilon(1)d(x_0)^{-\alpha} + Lo_\epsilon(1).
\end{split}
\end{equation*} 

Letting $\delta' \to 0$ we get rid of the term $\epsilon^{-2} O({\delta'}^{1 - \alpha})$. By~\eqref{propcorUV}, letting $\epsilon \to 0$ 
we handle the differential and integral terms by the smoothness of $\phi$ and the semicontinuity of $U$ and $V$ together with Fatou's Lemma. 
Finally, letting $h \to 0$ we conclude
$
\mathcal{E}_\delta^i(W, \phi, x_0) \leq 0,
$
which is the desired viscosity inequality for $W$. 

Now we deal with the case $x_0 \in \Gamma_{in}$. By definition of $U, V$ there exist sequences $(x_k)_k, (y_k)_k$ of points of $\Omega$ such that 
$x_k, y_k \to x_0$ and $U(x_k) \to U(x_0), V(y_k) \to V(x_0)$. Denote $\eta_k = \min \{ d(x_k), d(y_k)\}$, $\epsilon_k = |x_k - y_k|$ (which
can be taken strictly possitive for all $k$) and define the functions 
\begin{equation*}
\tilde{\phi}_k(x,y) = \phi(x) + \frac{|x - y|^2}{\epsilon_k} \quad \mbox{and} \quad \phi_k(x,y) 
= \tilde{\phi}_k(x,y) - \eta_k \zeta(x) - \eta_k \zeta(y).
\end{equation*}

Doubling variables, consider the penalization
\begin{equation*}
(x,y) \mapsto \Phi(x,y) := U(x) - V(y) - \phi_k(x, y), 
\end{equation*}
which is a function attaining its maximum in $(B_\delta(x_0) \cap \Omega) \times (B_\delta(x_0) \cap \Omega)$ at a point $(\bar{x}, \bar{y})$. Using 
the inequality $\Phi(\bar{x}, \bar{y}) \geq \Phi(x_k, y_k)$ we conclude
\begin{equation*}
|\bar{x} - \bar{y}|^2/\epsilon_k \to 0, \ 
\bar{x}, \bar{y} \to x_0, \ \ U(\bar{x}) \to U(x_0), \ V(\bar{y}) \to V(x_0) \quad \mbox{as} \ k \to \infty,  
\end{equation*}
and we are in the situation of an interior maximum point. As above, consider $0 < \delta' < c_k := \min \{d(\bar{x}), d(\bar{y})\}$ and consider the
respective viscosity inequalities for $U$ and $V$ on $\bar{x}$ and $\bar{y}$. 
At this point we use inequality~\eqref{testingpropUV2} and its corresponding version for supersolutions, and introducing the parameter 
$h > 0$ as in the previous case to avoid the sup terms we conclude
\begin{equation*}
\mathcal{E}_{\delta', \beta_h}(\bar{x}, U, \tilde{\phi}_k(\cdot, \bar{y})) \leq -\tilde{c}_0 d(\bar{x})^{-\alpha} \ \mbox{and} \
\mathcal{E}_{\delta', \beta_h}(\bar{y}, V, -\tilde{\phi}_k(\bar{x}, \cdot)) + h \geq \tilde{c}_0 d(\bar{y})^{-\alpha},
\end{equation*}
for some $\tilde{c}_0 > 0$ independent of $k$. Substracting both inequalities and using the definition of $\tilde{\phi}$, we can make 
similar computations as in~\eqref{driftcorpropUV}, ~\eqref{nonlocalcorpropUV} and~\eqref{delta'bydelta} relative to this case, to conclude
\begin{equation*}
\begin{split}
& - \mathcal{I}_{\Omega, \delta} [\phi](\bar{x}) - J_{int}^\delta + \mathcal{H}_i(\bar{x}, D\phi(\bar{x})) \\
\leq & \ h + (o_k(1) - \tilde{c}_0) c_k^{-\alpha} + \epsilon_k^{-2} o_{\delta'}(1) + L o_k(1) \\
& \ + \eta_k \Big{[} -\mathcal{I}_{\Omega, \delta}[\zeta](\bar{x}) + \mathcal{H}_s(\bar{x}, D\zeta(\bar{x}))
- \mathcal{I}_{\Omega, \delta}[\zeta](\bar{y}) + \mathcal{H}_s(\bar{y}, D\zeta(\bar{y})) \Big{]}
\end{split}
\end{equation*} 

Letting $\delta' \to 0$ we get rid of the term $\epsilon_k^{-2} o_{\delta'}(1)$. Then, using Corollary~\ref{corlog} we conclude that for all 
$k$ large the term in squared brackets in the right hand side of the last expression is nonpossitive. 
Letting $k \to \infty$ and then $h \to 0$ we conclude as in the previous case.
\qed


\section{The Cone Condition.}
\label{coneconditionsection}

In this section we provide a proof for the well-known ``cone condition'': the fact that we can 
approximate the value of a subsolution to problem~\eqref{eq}-\eqref{Dirichlet} at a point $x_0 \in \partial \Omega$ through a sequence of points  
lying in a cone contained in $\Omega$ whose vertex is $x_0$.
We remark this condition is also adressed in the 
probabilistic approach refered as \textsl{nontangential upper semicontinuity} (see~\cite{Katsoulakis}), but
here we follow closely the lines of the corresponding property proved in~\cite{Barles-Rouy}.

\begin{prop}\label{conecondition}
Let $u$ be a bounded viscosity subsolution of~\eqref{eq}-\eqref{Dirichlet} and $\tilde{u}$ defined in~\eqref{deftildeuv}. Then, for each 
$x_0 \in \Gamma \cup \Gamma_{in}$ there exists a constant $C > 0$ and a sequence $(x_k)_k$ of points of $\Omega$ such that, as $k \to \infty$
\begin{eqnarray}\label{conesequence}
\left \{ \begin{array}{l} x_k \to x_0, \\ \tilde{u}(x_k) \to \tilde{u}(x_0), \\ d(x_k) \geq C|x_k - x_0|. \end{array} \right .
\end{eqnarray}
\end{prop}

\noindent
{\bf Proof:} Consider $\bar{\beta}$ in~\eqref{Gammaexplanation} relative to $x_0$ and denote $b = b_{\bar{\beta}}$. Since we have 
$b(x_0) \cdot Dd(x_0) > 0$, we take $r > 0$ small enough such that $b(x) \cdot Dd(x) > 0$ for all $x \in \bar{\Omega} \cap \bar{B}_r(x_0)$.
After rotation and translation, we can assume $x_0 = 0$ and $Dd(x_0) = e_n$ with $e_n = (0,...,0,1)$, implying in particular that 
$b_n(0) > 0$. Finally, denote 
$H_+ = \{(x', x_n) \in \R^n : x_n > 0\}$ and $A = \bar{H}_+ \cap \bar{\Omega} \cap \bar{B}_r$.

Recalling $U$ defined as in~\eqref{U}, by Proposition~\ref{propU} we have this function satisfies the equation 
\begin{equation*}
-\mathcal{I}_\Omega[w] - b \cdot Dw \leq 0 \quad \mbox{on} \ A.
\end{equation*}

By a simple scaling argument, we conclude the function $y \mapsto U(\gamma y)$ defined in $\gamma^{-1} A$ satisfies the equation
\begin{equation}\label{transportgamma}
-\gamma^{1 - \alpha} \mathcal{I}_{\gamma^{-1}\Omega} [w](y) - b_\gamma (y) \cdot Dw(y) \leq 0 \quad \mbox{on} \ \gamma^{-1}A,
\end{equation}
where $b_\gamma(y) = b(\gamma y)$ for each $y \in \gamma^{-1}A$.
Thus, $\bar{w} : \bar{H}_+ \to \R$ defined as
\begin{eqnarray*}
\bar{w}(x) = \limsup \limits_{\gamma \to 0, z \to x} U (\gamma z)
\end{eqnarray*}
is a viscosity subsolution for the problem
\begin{eqnarray*}\label{transport}
- b_n(0) \frac{\partial w}{\partial y_n} - b'(0) \cdot D_{y'}w = 0 \quad \mbox{in} \ \bar{H}_+,
\end{eqnarray*}
by classical arguments in half-relaxed limits applied over the equation~\eqref{transportgamma}. It is worth remark that this equation holds up 
to the boundary and that $b_n(0) > 0$.

The maximal solution for the last transport equation with terminal data $\bar{w}(y', 1)$ (when we cast $y_n$ as the ``time'' variable) 
is given by the function
$$
W(y', y_n) = \bar{w}(y' - b_n(0)^{-1} b'(0)(y_n - 1), 1).
$$ 

Since $W$ is maximal, we have $\bar{w}(y) \leq W(y)$ when $0 \leq y_n \leq 1$. Now, by definition it is clear that $\bar{w}$ is upper semicontinuous
and then $\bar{w}(0) = U(0)$, meanwhile by the upper semicontinuity of $u$ at the boundary and the continuity of the distance function we 
have $\bar{w}(y) \leq U(0)$ for all $y \in H_+$. Then, recalling $U(0)=\tilde{u}(0)$, we conclude that
$$
\tilde{u}(0) = \bar{w}(0) \leq W(0) = \bar{w}(b_n(0)^{-1} b'(0), 1) \leq \tilde{u}(0),
$$
this is $\tilde{u}(0) = \bar{w}(x_b)$, with $x_b = (b_n(0)^{-1} b'(0), 1)$. By the very definition of $\bar{w}$, we have the existence
of sequences $\gamma_k \to 0$, $z_k \to x_b$ such that $x_k := \gamma_k z_k$ satisfies $x_k \to 0$ and $\tilde{u}(x_k) \to \tilde{u}(0)$. 

Note that by definition of the sequence $(x_k)_k$ we have 
$
x_k = \gamma_k x_b + o(\gamma_k).
$
Using this, we perform a Taylor expansion on $d(x_k)$, obtaining the existence of a point $\bar{x}_k \in H_+$ with $\bar{x}_k \to 0$ as 
$k \to \infty$ such that
\begin{equation*}
d(x_k) = Dd(\bar{x}_k) \cdot (\gamma_k x_b + o(\gamma_k)).
\end{equation*}

Hence, since $Dd(0) = e_n$ we conclude $d(x_k) = \gamma_k + o(\gamma_k)$. Thus, using the estimates for $x_k$ and $d(x_k)$ we get that
$
d(x_k) \geq (4|x_b|)^{-1} |x_k|,
$
for all $k$ large enough. Recalling that $x_0 = 0$, we conclude that $(x_k)_k$ is the sequence satisfying~\eqref{conesequence}.
\qed

Now we state the the analogous result for supersolutions, valid on $\Gamma_{in}$.

\begin{prop}\label{superconecondition}
Let $v$ be a bounded viscosity supersolution of~\eqref{eq}-\eqref{Dirichlet} and let $\tilde{v}$ as in~\eqref{deftildeuv}. Then, for each 
$x_0 \in \Gamma_{in}$, there exists a sequence $(x_k)_k$ of points of $\Omega$ satisfying~\eqref{conesequence} relative to $\tilde{v}$.
\end{prop}

\noindent
{\bf Proof:} Following the arguments in the previous proposition, we consider this time the function
\begin{equation*}
\underline{w}(x) = \liminf \limits_{\gamma \to 0, z \to x} V(\gamma z), 
\end{equation*}
where $V$ is defined in~\eqref{V}. Using Proposition~\ref{propUV} relative to $V$, we prove $\underline{w}$ is a viscosity supersolution to 
the transport equation
\begin{eqnarray}\label{eqtransportsuper}
\sup \limits_{\beta \in \mathcal{B}}  \{ -b_\beta(0) \cdot Dw \} = 0 \quad \mbox{in} \ \bar{H}_+.
\end{eqnarray}

Since $x_0 \in \Gamma_{in}$, $(b_\beta)_n(0) > 0$ for all $\beta \in \mathcal{B}$ and then $\tilde{b}_\beta(0) = b_\beta(0) / (b_\beta)_n(0)$
is well defined.
Now, if we denote $W$ a minimal solution to~\eqref{eqtransportsuper}, then it is easy to see that $W$ is a solution to
\begin{eqnarray*}
- \frac{\partial w}{\partial y_n} + \sup \limits_{\beta \in \mathcal{B}} \{ \tilde{b}_\beta'(0) \cdot D_{y'}w \} = 0 \quad \mbox{in} \ \bar{H}_+.
\end{eqnarray*}

By optimal control arguments (see~\cite{Bardi-Capuzzo},~\cite{Barles-book}), it is known the minimal solution of this equation with terminal 
data $W(y', 1) = \underline{w}(y',1)$ is
$$
W(y', y_n) = \inf \limits_{\beta(\cdot)} \{ \underline{w}(Y_{y'}((1 - y_n), \beta), 1)\}.
$$
where $t \mapsto Y_{y'}(t, \beta(\cdot))$ is the solution of the equation
$
\dot{Y}(t) = \tilde{b}'(0, \beta(t))
$
with initial condition $Y(0) = y'$. Note that in particular, we have 
\begin{eqnarray}\label{superconeconditionestimateI}
|Y(t)| \leq |y'| + t ||b_\beta(0)||_{L^\infty(\mathcal{B})}, 
\end{eqnarray}
for all $y' \in \R^{n - 1}$ and all $t \in [0,1]$. At this point, arguing in the same way as in the proof of 
Proposition~\ref{conecondition}, we conclude
$$
\tilde{v}(0) = \inf \limits_{\beta(\cdot)} \{ \underline{w}(Y_{0}(1, \beta), 1)\}.
$$

By this expression, we can take a minimizing sequence of controls $\beta_k(\cdot)$ and denoting $z_k = (Y_{0}(1, \beta_k), 1)$ we have 
\begin{equation*}
\underline{w}(z_k) \to \tilde{v}(0) \quad \mbox{as} \ k \to \infty.
\end{equation*}

Recalling that $V = \tilde{v} - Ad^{1 - \alpha}$, by the very definition of $\underline{w}$ and the continuity of the distance function, 
for each $k$ there exists sequences $\gamma^j_k \to 0$, $z_k^j \to z_k$ such that 
$\tilde{v}(\gamma_k^j z_k^j) \to \underline{w}(z_k)$ as $j \to \infty$. Note that by~\eqref{superconeconditionestimateI} we have the sequence
$(z_k^j)_{k,j}$ is bounded and therefeore, using a diagonal argument we conclude the existence of sequences $\bar{\gamma}_k \to 0$,
$\bar{z}_k \to \bar{z}$ with $\bar{z} = (\bar{z}', 1)$ and $||\bar{z}'|| \leq ||b_\beta(0)||_{L^\infty(\mathcal{B})}$ such that
$\tilde{v}(\bar{\gamma}_k \bar{z}_k) \to \tilde{v}(0)$. Arguing in a similar way as in the end of the previous proposition, we conclude
the sequence $(x_k)_k$ defined by $x_k = \bar{\gamma}_k \bar{z}_k$ is the desired sequence satisfying~\eqref{conesequence} relative to $\tilde{v}$.
\qed

Propositions~\ref{conecondition} and~\ref{superconecondition} are sufficient to conclude Theorems~\ref{teo1} and~\ref{teo2}. However, 
we provide an alternative proof for Case III in section~\ref{proofteo} below, which is more natural according the expected behavior
of the solutions near $\Gamma_{in}$.
For this, we present next a lemma in the flavour of the last two propositions but concerning the difference $\tilde{u} - \tilde{v}$.

\begin{lema}\label{lemaw}
Let $u, v$  be respectively bounded viscosity sub and supersolution to problem~\eqref{eq}-\eqref{Dirichlet}, $\tilde{u}, \tilde{v}$ as 
in~\eqref{deftildeuv} and define $\omega = \tilde{u} - \tilde{v}$.
Let $r > 0$ and $x_0 \in \Gamma_{in}$ such that $\omega(x_0) \geq \omega(x)$ for all $x \in B_r(x_0) \cap \Gamma_{in}$. Then, there exists a
sequence $(x_k)_k$ of points of $\Omega$ such that $\omega(x_k) \to \omega(x_0)$. 
\end{lema}
 
\noindent
{\bf Proof:} By contradiction we assume there exists $r' > 0$ such that $\omega(x_0) > \omega(x)$ for all 
$x \in B_{r'}(x_0) \cap \Omega$, and we may assume that $r' \leq r$. Let $U$ as in ~\eqref{U}, $V$ as in~\eqref{V} and 
$W = U - V$. Since this function is upper semicontinuous, by the contradiction assumption there exists $m > 0$ such that $W(x_0) \geq W(x) + m$ 
for all $x \in B_r(x_0) \cap \Omega$, taking $r$ smaller if it is necessary. By the hypothesis over $\omega$ on the boundary, for each $\epsilon > 0$, 
$x_0$ is the maximum point of the function
$
x \mapsto W(x) + \epsilon^{-1} d(x)
$
in $B_{m\epsilon}(x_0) \cap \bar{\Omega}$, and then we can use this as a testing for $W$ at $x_0$ by Proposition~\ref{propW}, concluding 
for each $0 < \delta < m\epsilon$ that
\begin{equation}\label{testinglemaw}
\mathcal{E}_\delta^i(W, -\epsilon^{-1}d, x_0) \leq 0. 
\end{equation}

However, note that there exists $c_0 > 0$ such that $b_\beta(x_0) \cdot Dd(x_0) \geq c_0$ for 
all $\beta \in \mathcal{B}$ and this implies that
\begin{equation}\label{ineqlemaw1}
\mathcal{H}_i(x_0, -\epsilon^{-1}Dd(x_0)) \geq \epsilon^{-1}c_0.
\end{equation}

On the other hand, since $W(x) \leq W(x_0)$ for $x \in B_r(x_0)$ we have
\begin{equation*}\label{ineqlemaw2}
\mathcal{I}_\Omega^\delta[W](x_0) \leq \mathcal{I}_\Omega^r[W](x_0) \leq 2||W||_{L^\infty(\bar{\Omega})} r^{-\alpha}.
\end{equation*}

Finally, using this inequality and~\eqref{ineqlemaw1} in~\eqref{testinglemaw}, we conclude
\begin{equation*}
- \epsilon^{-1} o_\delta(1) - 2||W||_{L^\infty(\bar{\Omega})} r^{-\alpha} + c_0 \epsilon^{-1} \leq 0,
\end{equation*}
which is a contradiction after doing $\delta \to 0$ and then $\epsilon \to 0$.
\qed


\section{Proof of the Main Results.}
\label{proofteo}

\subsection{Proof of Theorem~\ref{teo1}.}
Our arguments rely on the redefined functions $\tilde{u}, \tilde{v}$, but we omit superscript ``$\sim$'' for simplicity. By contradiction, we assume
\begin{equation*}
M :=\sup \{ u(x) - v(x) : x \in \bar{\Omega} \}> 0. 
\end{equation*}

This supremum is attained at some point $x_0 \in \bar{\Omega}$ by the upper semicontinuity of $u - v$ in $\bar{\Omega}$ and recalling 
Corollary~\ref{corlemaGammaout} we have $x_0 \notin \Gamma_{out}$. From now, we proceed classically doubling variables and penalizing the difference 
$u(x) - u(y)$ by suitable functions $\phi(x,y)$, distinguishing where is the maximum point.

\medskip

\noindent
{\bf Case I: There is a maximum point $\Omega.$} 

\medskip

For $\gamma, \eta > 0$ denote $\phi_{\gamma, \eta}(x,y) = \frac{|x - y|^2}{\eta^2} + \gamma|x - x_0|^2 $ and penalize as
\begin{eqnarray*}\label{penCaseI}
(x,y) \mapsto \Phi_{\gamma, \eta}(x,y) := u(x) - v(y) - \phi_{\gamma, \eta}(x,y).
\end{eqnarray*}

For all $\eta, \gamma$ suitable small, $\Phi_{\gamma, \eta}$ attains its maximum point over $\bar{\Omega} \times \bar{\Omega}$ at 
$(\bar{x}_{\gamma, \eta}, \bar{y}_{\gamma, \eta})$. Here and in the next cases, we drop the dependence of the maximum points and test 
functions on the introduced parameters to avoid overwritten expressions. 
By the inequality $\Phi(\bar{x}, \bar{y}) \geq \Phi(x_0, x_0)$ we conclude
\begin{eqnarray*}\label{ineqCaseI}
0 < M \leq u(\bar{x})- v(\bar{y}) - \frac{|\bar{x} - \bar{y}|^2}{\eta^2} - \gamma|\bar{x} - x_0|^2
\end{eqnarray*}
and from this, it is easy to conclude, as $\eta \to 0$
\begin{equation}\label{propertiesCaseI}
|\bar{x} - \bar{y}|^2/\eta^2 \to 0; \quad \bar{x}, \bar{y} \to x_0; \quad u(\bar{x}) - v(\bar{y}) \to M, 
\end{equation}

In particular, if $\gamma >  0$ is fixed, $d(\bar{x}), d(\bar{y}) \geq d(x_0)/2 > 0$
for all $\eta$ suitably small. Using this penalization as test functions for $u$ and $v$, we conclude
$
E_\delta(u^\varphi, \phi(\cdot, \bar{y}), \bar{x}) \leq 0 \leq E_\delta(v_\varphi, -\phi(\bar{x}, \cdot) , \bar{y})
$
for all $0 < \delta < d(x_0)/2$, where we have used notation~\eqref{Enonlinear}. By definition of $H$ and notation~\eqref{Eevaluation},
for each $h > 0$ there exists $\beta_h \in \mathcal{B}$ such that
\begin{eqnarray*}
E_{\beta_h, \delta}(u^\varphi, \phi(\cdot, \bar{y}), \bar{x}) \leq 0 \quad \mbox{and} \quad 
-h \leq E_{\beta_h, \delta}(v_\varphi, -\phi(\bar{x}, \cdot) , \bar{y}).  
\end{eqnarray*}

Then, substracting these inequalities, we conclude
\begin{equation}\label{testingCaseI}
\begin{split}
& \lambda (u(\bar{x}) - v(\bar{y})) - h \\
\leq & \ f_{\beta_h}(\bar{x}) - f_{\beta_h}(\bar{y}) + 2 \gamma b_{\beta_h}(\bar{x}) \cdot (\bar{x} - x_0)
+ 2\eta^{-2} (b_{\beta_h}(\bar{x}) - b_{\beta_h}(\bar{y})) \cdot (\bar{x} - \bar{y}) \\
& + \mathcal{I}_\delta[\phi(\cdot, \bar{y})](\bar{x}) + \mathcal{I}_\delta[\phi(\bar{x}, \cdot)](\bar{y}) 
+I_{int}^\delta + I_{ext}^\delta + I_{int, \bar{y}}^\delta + I_{int, \bar{x}}^\delta,
\end{split}
\end{equation}
where, using notation~\eqref{defD's} and for $\delta' > 0$ we have denoted
\begin{equation}\label{I's}
\begin{split}
& I_{int}^{\delta'} = \int \limits_{D_{int} \setminus B_{\delta'}} [u(\bar{x} + z) - v(\bar{y} +z) - (u(\bar{x}) - v(\bar{y}))]K^\alpha(z)dz, \\
& I_{ext}^{\delta'} = 
\int \limits_{D_{ext} \setminus B_{\delta'}} [\varphi(\bar{x} + z) - \varphi(\bar{y} +z) - (u(\bar{x}) - v(\bar{y}))]K^\alpha(z)dz, \\  
& I_{int, \bar{y}}^{\delta'} = 
\int \limits_{D_{int}^{\bar{y}} \setminus B_{\delta'}} [\varphi(\bar{x} + z) - v(\bar{y} +z) - (u(\bar{x}) - v(\bar{y}))]K^\alpha(z)dz,\\
& I_{int, \bar{x}}^{\delta'} = 
\int \limits_{D_{int}^{\bar{x}} \setminus B_{\delta'}} [u(\bar{x} + z) - \varphi(\bar{y} +z) - (u(\bar{x}) - v(\bar{y}))]K^\alpha(z)dz.
\end{split}
\end{equation}

We need to estimate each term in the right-hand side of~\eqref{testingCaseI}. First, by the continuity of $f$, condition~\eqref{L} and the properties
in~\eqref{propertiesCaseI}, we conclude
\begin{equation}\label{estimateCaseI1}
f_{\beta_h}(\bar{x}) - f_{\beta_h}(\bar{y}) + 2 \gamma b_{\beta_h}(\bar{x}) \cdot (\bar{x} - x_0)
+ 2\eta^{-2} (b_{\beta_h}(\bar{x}) - b_{\beta_h}(\bar{y})) \cdot (\bar{x} - \bar{y}) = o_\eta(1). 
\end{equation}

Now we take care about the nonlocal terms. By definition of $\phi$ we have
\begin{equation}\label{estimateCaseI2}
\mathcal{I}_\delta[\phi(\cdot, \bar{y})](\bar{x}) + \mathcal{I}_\delta[\phi(\bar{x}, \cdot)](\bar{y}) = (\eta^{-2} + \gamma)O(\delta^{1 - \alpha}). 
\end{equation}

Note the sets $D_{int}^{\bar{x}}, D_{int}^{\bar{y}}$ vanish as $\eta \to 0$ and are away from the origin 
at least at distance $d(x_0)/2$ for all $\eta$ suitable small. This implies 
\begin{equation}\label{estimateCaseI3}
I_{int, \bar{x}}^\delta, I_{int, \bar{y}}^\delta = o_\eta(1) d(x_0)^{-\alpha}
\end{equation}
by the first fact in~\eqref{propertiesCaseI}. Besides, $D_{ext}$ is also away from the origin and for each $z \in D_{ext}$, using 
again~\eqref{propertiesCaseI} and the continuity of $\varphi$ we have
\begin{equation}\label{to-M}
\varphi(\bar{x} + z) - \varphi(\bar{y} + z) - (u(\bar{x}) - v(\bar{y})) \to -M
\end{equation}
as $\eta \to 0$. Hence, applying Dominated Convergence Theorem we conclude
\begin{equation}\label{estimateCaseI4}
I_{ext}^\delta \to -M \int_{(\Omega - x_0)^c} K_\alpha(z)dz \quad \mbox{as} \ \delta, \eta \to 0.
\end{equation}

Finally, using $(\bar{x}, \bar{y})$ is a maximum point of $\Phi$ in $\bar{\Omega} \times \bar{\Omega}$ it is easy to see that
for each $z \in D_{int}$ we have the inequality
\begin{equation*}
u(\bar{x} + z) - v(\bar{y} + z) - (u(\bar{x}) - v(\bar{y})) \leq \gamma (|\bar{x} + z - x_0|^2 - |\bar{x} - x_0|^2),
\end{equation*}
which allows to conclude
\begin{equation}\label{estimateCaseI5}
I_{int}^\delta = O(\gamma). 
\end{equation}

Applying this last estimate together with~\eqref{estimateCaseI1},~\eqref{estimateCaseI2},~\eqref{estimateCaseI3},~\eqref{estimateCaseI4}  
in~\eqref{testingCaseI} we get
\begin{equation*}
\begin{split}
& \lambda (u(\bar{x}) - v(\bar{y})) + M \int_{(\Omega - x_0)^c} K_\alpha(z)dz \\
\leq & \ h + o_\delta(1) + o_\eta(1) + (\eta^{-2} + \gamma)O(\delta^{1 - \alpha}) + o_\eta(1) d(x_0)^{-\alpha} + O(\gamma).  
\end{split}
\end{equation*}

Hence, using~\eqref{propertiesCaseI} and let $\delta \to 0$, $\eta \to 0$, $\gamma \to 0$ and $h \to 0$ we arrive to
\begin{equation*}
M \Big{(} \lambda + \int_{(\Omega - x_0)^c} K^\alpha(z)dz \Big{)} \leq 0, 
\end{equation*}
which is a contradiction because of assumption~\eqref{measure}.

\medskip

\noindent
{\bf Case II: There is a maximum point in $\Gamma$.} 

\medskip

Note that by Proposition~\ref{lemaGammaout}, $u(x_0) \leq \varphi(x_0)$ and by the 
contradiction assumption, the only possibility is $v(x_0) < \varphi(x_0)$. This implies that $v_\varphi(x_0) = v(x_0)$. 

Consider the sequence $(x_k)_k$ satisfying~\eqref{conesequence} relative to $u$ in Theorem~\ref{conecondition} and define 
$\eta_k = |x_k - x_0|$ and 
$\nu_k = (x_k - x_0)/\eta_k$. By the cone condition (Proposition~\ref{conecondition}), we can assume the sequence in such a way $\nu_k \to \nu_0$ and
$Dd(x_0) \cdot \nu_0 > 0$. Following the arguments of Theorem 7.9 in~\cite{usersguide}, we use the function
$$
\phi(x,y) = |\eta_k^{-1} (x-y) - \nu_0|^2 + \gamma|y - x_0|^2,
$$
where $\gamma > 0$, and then we write the penalization as 
$$
(x,y) \mapsto \Phi(x,y) := u(x) - v_\varphi(y) - \phi(x,y).
$$

With this, we have a maximum point $(\bar{x},\bar{y}) \in \bar{\Omega} \times \R^n$ for $\Phi$. Using the inequality 
$\Phi(\bar{x}, \bar{y}) \geq \Phi(x_k, x_0)$ we get
\begin{equation}\label{ineqCaseII}
\begin{split}
\gamma |\bar{y} - x_0|^2 + |\eta_k^{-1} (\bar{x}- \bar{y}) - \nu_0|^2 
\leq u(\bar{x}) - v_\varphi(\bar{y}) - (u(x_k) - v(x_0)) +  |\nu_k - \nu_0|^2. 
\end{split}
\end{equation}

Since $u, v$ are bounded, $|\eta_k^{-1}(\bar{x}- \bar{y}) - \nu_0|^2$ remains bounded too and then 
$|\bar{x} - \bar{y}| \to 0$ as $k \to \infty$. Using the semicontinuity of $u$ and $v_\varphi$ and since $u(x_k) \to u(x_0)$, 
we obtain from~\eqref{ineqCaseII}, up to a subsequence 
\begin{equation}\label{propertiesCaseII}
\bar{x}, \bar{y} \to x_0,  \ u(\bar{x}) \to u(x_0), \ v_\varphi(\bar{y}) \to v(x_0), \ \phi(\bar{x}, \bar{y}) \to 0; \quad \mbox{as} \  
k \to \infty.
\end{equation}

With this we have $\bar{x}, \bar{y}$ are valid test points for $u$ and $v$, respectively, using the penalizations introduced as test functions. 
In fact, we have $\bar{y} \in \bar{\Omega}$, otherwise we get 
$v_\varphi(\bar{y}) \to \varphi(x_0) > v(x_0)$ which is a contradiction with the third statetment in~\eqref{propertiesCaseII}. 
Since $\bar{y}$ is a global minimum point of
$
y \mapsto v_\varphi(y) - (u(\bar{x}) - \phi(\bar{x}, y)),
$
we can use $\bar{y}$ as a test point for $v_\varphi$ even if it is on the boundary because if this is the case, $v_\varphi(\bar{y}) < \varphi(\bar{y})$
by the continuity of $\varphi$. On the other hand, again by~\eqref{ineqCaseII} we have
\begin{equation}\label{orderofbarx-bary}
\bar{x} = \bar{y} + \eta_k (\nu_0 + o_k(1)). 
\end{equation}

By a Taylor expansion of $d(\bar{x})$ we conclude the existence of $\tilde{y} \to x_0$ as $k \to \infty$ such that 
$
d(\bar{x}) \geq d(\bar{y}) + \eta_k(Dd(\tilde{y}) \cdot \nu_0 + o_k(1)).
$
Since $Dd(x_0) \cdot \nu_0 > 0$ and by the continuity of $Dd$ we conclude $d(\bar{x}) > 0$ for all $k$ large enough, meaning
$\bar{x} \in \Omega$ and then we use it as a test point for $u$. 
Then, since $u$ and $v$ are respective viscosity sub and supersolutions to problem~\eqref{eq}-\eqref{Dirichlet}, for all $0 < \delta < d(\bar{x})$
we can write
$
E_\delta(u, \phi(\cdot, \bar{y}), \bar{x}) \leq 0 \leq E_\delta(v, -\phi(\bar{x}, \cdot), \bar{y})
$
and then, for all $h > 0$ there exists $\beta_h \in \mathcal{B}$ such that
\begin{eqnarray*}
E_{\delta, \beta_h}(u, \phi(\cdot, \bar{y}), \bar{x}) \leq 0 \quad \mbox{and} \quad E_{\delta, \beta_h}(v, -\phi(\bar{x}, \cdot), \bar{y}) \geq -h.
\end{eqnarray*}

Substracting these last inequalities we arrive to an inequality in the same way as we got~\eqref{testingCaseI}.
Using \eqref{propertiesCaseII} and \eqref{L} we argue as in the previous case, concluding this time the inequality
\begin{equation}\label{testingCaseII}
\lambda M \leq \ h + \eta_k^{-4} O(\delta^{1 - \alpha}) + o_k(1) + O(\gamma) + 
I_{int}^\delta + I_{ext}^\delta + I_{int, \bar{y}}^\delta + I_{int, \bar{x}}^\delta,
\end{equation}
where we have used notation~\eqref{I's}.

In what follows, we estimate the nonlocal terms in the right-hand side of the last inequality. 
For all $r_1 > 0$ small enough, since $u(x_0) \leq \varphi(x_0)$ and the upper semicontinuity of $u$ we have
\begin{equation*}
u(\bar{x} + z) - \varphi(\bar{y} +z) - (u(\bar{x}) - v(\bar{y})) \leq 0, \quad \forall \ z \in B_{r_1} \cap D_{int}^{\bar{x}},
\end{equation*}
and then $I_{int, \bar{x}}^\delta \leq I_{int, \bar{x}}^{r_1}$. This estimate allows us to get rid of the dependence of $\delta$ of 
$I_{int, \bar{x}}^\delta$, and since $D_{int}^{\bar{x}}$ vanishes as $k \to \infty$ we conclude
\begin{equation}\label{estimateCaseII1}
I_{int, \bar{x}}^\delta = o_k(1).
\end{equation}

By~\eqref{propertiesCaseII} and the continuity of $\varphi$ we have
\begin{equation*}
\varphi(\bar{x} + z) - \varphi(\bar{y} +z) - (u(\bar{x}) - v(\bar{y})) \leq 0, \quad \forall \ z \in B_{r_1} \cap D_{ext}, 
\end{equation*}
and then $I_{ext}^\delta \leq I_{ext}^{r_1}$. Using~\eqref{to-M} over $I_{ext}^{r_1}$ we conclude
\begin{equation}\label{estimateCaseII2}
I_{ext}^{\delta} \leq -M \int_{(\Omega - x_0)^c \setminus B_{r_1}} K^\alpha(z)dz + o_k(1).
\end{equation}

At this point, by~\eqref{orderofbarx-bary} we claim the set $D_{int}^{\bar{y}}$ is away from the origin uniformly in $k$, 
postponing the proof of this claim until the end of this Case. Using this and the fact that $D_{int}^{\bar{y}}$ vanishes as $k \to \infty$, we 
conclude that
\begin{equation}\label{estimateCaseII3}
I_{int, \bar{y}}^\delta = o_k(1).  
\end{equation}

Finally, using that $(\bar{x}, \bar{y})$ is a maximum point for $\Phi$, in a similiar way 
as in~\eqref{estimateCaseI5} we get $I_{int}^\delta = O(\gamma)$ . Using this last estimate 
together with \eqref{estimateCaseII1}, \eqref{estimateCaseII2}, \eqref{estimateCaseII3} in~\eqref{testingCaseII} we arrive to
\begin{equation*}
M \Big{(} \lambda + \int_{(\Omega - x_0)^c \setminus B_{r_1}} K^\alpha(z)dz \Big{)} \leq \ h + \eta_k^{-4} O(\delta^{1 - \alpha}) + O(\gamma) + o_k(1).
\end{equation*}

From this, letting $\delta \to 0$, $k \to \infty$, $\gamma \to 0$ and $h \to 0$ we arrive to a contradiction with~\eqref{measure}, taking $r_1$ 
smaller if it is necessary.

Now we address the claim leading to~\eqref{estimateCaseII3}. Assume that there exists a sequence $z_k \in D_{int}^{\bar{y}}$ such that $z_k \to 0$. 
By definition, there exists
$a_k \in \Omega$ and $b_k \in \Omega^c$ such that $z_k = a_k - \bar{y} = b_k - \bar{x}$ and by the first statement in~\eqref{propertiesCaseII} we
have $a_k, b_k \to x_0$. Now, applying~\eqref{orderofbarx-bary} we conclude $b_k = a_k + \eta_k^2(\nu_0 + o_k(1))$. Taking $k$ large we conclude
$b_k \in \Omega$, which is a contradiction.

\medskip
\noindent
{\bf Case III: There is a maximum point in $\Gamma_{in}$.}

\medskip

By Propositions~\ref{conecondition} and~\ref{superconecondition} we have cone condition both for sub and supersolutions in $\Gamma_{in}$ and
then we could argue exactly as in the previous case assuming $v(x_0) < u(x_0) \leq \varphi(x_0)$ or $v(x_0) < \varphi(x_0) \leq u(x_0)$, and
just reversing the roles of $u$ and $v$ in the case $\varphi(x_0) \leq v(x_0) < u(x_0)$. However, we provide here an alternative proof 
without assuming any predefined order among $u(x_0), v(x_0)$ and $\varphi(x_0)$. 

If we denote $\omega = u - v$, we claim that this function is a viscosity subsolution to 
\begin{equation}\label{eqCaseIII}
\lambda w - \mathcal{I}_\Omega[w] + \mathcal{H}_i(x, Dw) = 0  
\end{equation}
at $x_0$. Thus, we can use the zero constant function as a test function 
for $\omega$ at $x_0$ to get
\begin{equation*}
\lambda \omega(x_0) - \mathcal{I}_\Omega^\delta [\omega](x_0) \leq 0,
\end{equation*}
and since $\lambda \geq 0$ we conclude $- \mathcal{I}_\Omega^\delta [\omega](x_0) \leq 0$. However, since $x_0$ is maximum for $\omega$ 
in $\bar{\Omega}$ we conclude the reverse inequality and then $\mathcal{I}_{\Omega}^\delta[\omega](x_0) = 0$. This implies that there exists
$x_0' \in \Omega$ such that $\omega(x_0') = M$ and then we fall into Case I, concluding the contradiction and the proof of the Theorem.

Now we give a brief proof of the claim since it relies in arguments we already used in section~\ref{propertiessection}.
First, it is easy to prove that $\omega$ is a viscosity subsolution to~\eqref{eqCaseIII} in the set $A = \Sigma_a \cap \{ x : \omega(x) \geq 0 \}$.
The proof of this fact follows the same lines of the first case presented in Proposition~\ref{propW}, doubling variables and penalizing with 
the function $|x - y|^2 / \epsilon^2$ as in~\eqref{penalizationcorpropUV}. Now, since $\omega(x_0)$ is maximum of $\omega$ in $\bar{\Omega}$,  
Lemma~\ref{lemaw} holds, concluding the existence of a sequence $(x_k)_k$ of points in $\Omega$ such that $\omega(x_k) \to \omega(x_0)$. 
Therefore, since $\omega(x_0) = M > 0$, we have $\omega(x_k) = M + o_k(1)$ for all $k$ large. Then, if $\phi$ is a smooth function such that $x_0$ is 
a maximum point of $\omega -\phi$, we can penalize this testing with the function $\eta_k \zeta$, where $\eta_k = d(x_k)$. This penalization pushes
the testing point in $A$ and then we conclude the result using the equation~\eqref{eqCaseIII} valid in $A$ and letting $k \to \infty$
by using Corollary~\ref{corlog}.
\qed

\begin{remark}
In the proof of Case III presented above, the possibility of testing on $\Gamma_{in}$ no matter the boundary condition is, says that the points 
on $\Gamma_{in}$ ``behave'' like interior points.
\end{remark}

The following property is a direct consequence of Theorem~\ref{teo1}.

\begin{cor}
Assume~\eqref{ellipticity},~\eqref{measure},~\eqref{H},~\eqref{L} hold. Consider $\varphi_1, \varphi_2: \Omega^c \to \R$ continuous and 
bounded functions. Let $u$ be a bounded viscosity subsolution to problem~\eqref{eq}-\eqref{Dirichlet} with exterior data $\varphi=\varphi_1$
and $v$ a bounded viscosity supersolution to problem~\eqref{eq}-\eqref{Dirichlet} with exterior data $\varphi = \varphi_2$. 
If $\varphi_1 \leq \varphi_2$, then $u \leq v$ in $\Omega$.
Moreover, using Definition~\ref{deftildeuv}, we have $\tilde{u} \leq \tilde{v}$ in $\bar{\Omega}$.
\end{cor}


\subsection{Proof of Theorem~\ref{teo2}.}

By the smoothness of the boundary and assumption~\eqref{L}, we can consider a Lipschitz extension to all $\R^n$ of the function $b$ and 
continuous extension to $\R^n$ for $f$. With this, the corresponding extended Hamiltonian (which we still denote by $H$) can be understood as a 
function $H: \R^n \times \R^n \to \R$. Consider for each $\epsilon > 0$ continuous functions 
$\psi_{+}^\epsilon, \psi_{-}^\epsilon : \R^n \to \R$ with 
$\psi_+^\epsilon \geq \psi_-^\epsilon$ in $\R^n$ and $\psi_{\pm}^\epsilon = \varphi$ in $\Omega^c$. With this, for $w, l \in \R$ and 
$x, p \in \R^n$ we define the following Hamiltonian 
$$
H_\epsilon(x, w, p, l) = \min \{ (w - \psi_-^\epsilon(x)), \max \{ (w - \psi_+^\epsilon(x)), \lambda w - l + H(x, p) \} \},
$$
and its asociated problem
\begin{eqnarray}\label{eqauxteo2}
H_\epsilon(x, w(x), Dw(x), \mathcal{I}[w](x)) = 0, \ x \in \R^n.
\end{eqnarray}

It is easy to see that this problem admits comparison principle for bounded sub and supersolutions defined in all $\R^n$, see~\cite{Barles-Imbert}.
Recalling $\mu_0$ as in~\eqref{ellipticity}, let $R = ||\varphi||_\infty + \mu_0^{-1}||f||_\infty$ and consider the function 
$$
g(x) = 2R \mathbf{1}_{\bar{\Omega}}(x) + R \mathbf{1}_{\bar{\Omega}^c} (x).
$$

Note that $g \geq \varphi$ on $\Omega^c$ and that for each $x \in \Omega$ we have
\begin{equation*}
\lambda g(x) - \mathcal{I}_\Omega[g](x) - \int_{\Omega - x} [\varphi(x + z) - g(x)]K^\alpha(z)dz + H(x, 0) \geq 0,
\end{equation*}
concluding that $g$ is a supersolution to problem~\eqref{eqauxteo2}. In the same way can be proved that $-g$ is a subsolution for this problem.
Thus, Perron's Method applies, concluding the existence of a solution $w_\epsilon$ for this problem. 
Such solution satisfies $|w_\epsilon| \leq R$, which is an estimate independent of $\epsilon$. Now, considering $\psi_{\pm}^\epsilon$ such 
that $\psi_{\pm}^\epsilon \to \pm \infty$ in $\Omega$ as $\epsilon \to 0$, by the uniform boundeness of the function $w_\epsilon$, for all $x \in \R^n$ 
the quantities
$$
\bar{w}(x) := \limsup \limits_{\epsilon \to 0, y \to x} w_\epsilon (y), \quad 
\underline{w}(x) := \liminf \limits_{\epsilon \to 0, y \to x} w_\epsilon (y)
$$
are well defined. The half-relaxed limits method implies $\bar{w}, \underline{w}$ are bounded sub and supersolution 
to~\eqref{eq} in the sense of Definition~\ref{defsolproblem}, respectively. By the very definition of these functions, 
$\underline{w} \leq \bar{w}$ in $\R^n$. On the other hand, considering $u = \tilde{\bar{w}}, v = \tilde{\underline{w}}$ as in~\eqref{deftildeuv}
we apply Theorem~\ref{teo1} to get $u \leq v$ 
in $\bar{\Omega}$, concluding $u = v$ in $\R^n \setminus (\Gamma_{in} \cup \Gamma)$. However, for each $x_0 \in \Gamma \cup \Gamma_{in}$ we have
\begin{equation*}
v(x_0) \leq \liminf \limits_{x \in \Omega, x \to x_0} v(x) \leq \limsup \limits_{x \in \Omega, x \to x_0} v(x) 
= \limsup \limits_{x \in \Omega, x \to x_0} u(x) =  u(x_0)
\end{equation*}

It is easy to see $u$ and $v$ are respective viscosity sub and
supersolution to~\eqref{eq}-\eqref{Dirichlet} using Corollary~\ref{corlemaGammaout} and the fact that $u$ and $v$ differ from $\bar{w}$ and 
$\underline{w}$ in a set of zero Lebesgue measure. This concludes the function $u$ (equals $v$) is a $C(\bar{\Omega})$ viscosity solution, with 
$u = \varphi$ in $\bar{\Omega}^c \cup \Gamma_{out}$. The uniqueness follows directly by comparison.
\qed


\section{Remarks and Comments.}
\label{remarkssection}

Comparison principle for a class of fully nonlinear equations can be obtained by natural adaptation of the arguments presented in this paper.
For example, consider $0 < a < 1/2$ and $(\alpha_\beta)_{\beta \in \mathcal{B}}$ such that $\alpha_\beta \in (a, 1-a)$ for each $\beta \in \mathcal{B}$.
With this, we consider the nonlocal operator
\begin{equation}\label{Ibeta}
\mathcal{I}_\beta[\phi](x) = \int_{\R^n} \delta(\phi, x, z)K_\beta(z)|z|^{-(n + \alpha_\beta)}dz,  
\end{equation}
with $K_\beta: \R^n \to \R_+$ measurable function for all $\beta \in \mathcal{B}$ satisfying condition~\eqref{ellipticity} for constants 
$\Lambda, c_1$ 
independent of $\beta$. We also consider a continuous function $\lambda : \bar{\Omega} \times \mathcal{B} \to \R_+$ and introduce the following 
generalization of assumption~\eqref{measure}: For all $\beta \in \mathcal{B}$ 
\begin{equation}\tag{{\bf M'}}\label{M'}
\inf \limits_{x \in \bar{\Omega}} \Big{\{} \lambda_\beta(x) + \int \limits_{\Omega - x} K_\beta(z)dz \Big{\}} > 0.
\end{equation}

Let $\mathcal{C}$ be the set of real valued functions $\phi$ defined in all $\R^n$ for which~\eqref{Ibeta} is well defined for all 
$\beta \in \mathcal{B}$ and define $F : \R^n \times \R \times \mathcal{C} \times \R^n \to \R$ as
\begin{equation*}
F(x, z, \phi, p) = \sup \limits_{\beta \in \mathcal{B}} 
\{ \lambda_\beta(x)z -\mathcal{I}_\beta[\phi](x) -b_\beta(x) \cdot p - f_\beta(x)\}.
\end{equation*}

The fully nonlinear version of problem~\eqref{eq}-\eqref{Dirichlet} reads as
\begin{equation*}
\left \{ \begin{array}{rll} F(x, u(x), u, Du(x)) &= 0 \quad &\mbox{in} \ \Omega. \\
 u &= \varphi \quad &\mbox{in} \ \Omega^c. \end{array} \right .
\end{equation*}

Note that the definition of $\Gamma_{in}, \Gamma_{out}$ and $\Gamma$ depends only on the configuration of the drift term, so we shall use them 
exactly in the same form. By the assumption over the numbers $\alpha_\beta$, the order of the nonlocal operator is less than the order of the drift
so the results concerning section~\ref{technicalsection} can be readily adapted. On the other hand, the fractional nature of the operator is still
present in this case by the nonintegrability of the kernels, so results in section~\ref{propertiessection} and in paticular we can conclude 
cone condition for this type of equations. The conclusions in the contradiction arguments used to prove Theorem~\ref{teo1} can be obtained
in the same way in this setting by~\eqref{M'}.

The situation is different for nonlocal operators which are not uniformly elliptic. We address here two examples: For the first one, consider
$J \in L^1(\R^n)$ a positive function and denote
\begin{equation*}
\mathcal{J}[u](x) = \int \limits_{\R^n} \delta(u, x, z) J(z)dz.
\end{equation*}

The second one is the case of operators with censored jumps as
\begin{equation*}
\mathcal{I}_{+}[u](x) = \int \limits_{H_+} \delta(u, x, z) |z|^{-(n + \alpha)} dz,
\end{equation*}
with $\alpha \in (0,1)$ and $H_+ = \{(x', x_n) : x_n > 0\}$. It is known that in Dirichlet problems where $\mathcal{J}$ plays the diffusive role 
there exists loss of the boundary condition 
even in absence of drift term (see~\cite{Chasseigne}). A similar feature arises in equations related to operator $\mathcal{I}_+$, where despite the 
nonintegrability of the kernel, its null diffusive feature in some particular directions may create loss of the boundary condition at points where the 
censored direction is normal to the boundary. Hence, in both cases we cannot prescribe the value of the solution of the equation associated to these 
operators since definitions of $\Gamma_{out}$ and $\Gamma$ do not provide any information. However, we can obtain comparison results for these 
type of operators making stronger assumptions over $\Gamma_{out}$ and $\Gamma$. Define
\begin{equation*}
\Gamma_{out}' = \{ x \in \partial \Omega \ : \ \forall \ \beta \in \mathcal{B}, \ b_\beta(x) \cdot Dd(x) < 0 \}, 
\end{equation*}
consider the set $\Gamma_{in}$ as in the introduction and $\Gamma = \partial \Omega \setminus (\Gamma_{out}' \cup \Gamma_{in})$. 
Assuming the conditions

\noindent
{\bf (H')} $\Gamma_{out}'$, $\Gamma_{in}$ and $\Gamma$ are connected components of $\partial \Omega$, and

\medskip

\noindent
{\bf (H$\Gamma$)} For each $x \in \Gamma$, there exists $\bar{\beta}, \underline{\beta} \in \mathcal{B}$ such that
\begin{equation*}
b_{\underline{\beta}}(x) \cdot Dd(x) < 0 < b_{\bar{\beta}}(x) \cdot Dd(x),
\end{equation*}
then we can obtain strong comparison principle for Dirichlet problems associated to degenerate elliptic nonlocal operators. This is due to the 
stronger assumption over the controls on $\Gamma_{out}$ and $\Gamma$, since at one hand this implies Corollary~\ref{corlemaGammaout} holds, and 
on the other, by the fact we are still in the framework of low diffusive influence of the operator compared with the drift, technical results of 
section~\ref{propertiessection} and cone condition can be proved as they were presented here without substantial changes.



\section{Appendix.}
\label{app0}

\subsection{Proof of Lemma~\ref{lemaIlog}.} 
Let us first consider the illustrative case of a flat boundary, namely $\Omega = H_+ = \{(x',x_n) \in \R^{n - 1} \times \R : x_n > 0\}$. 
In this case, $d(x) = x_n$ and then for $x \in H_+$ we have
\begin{equation*}
\begin{split}
& -\mathcal{I}_{\Omega}[\zeta](x) = -\int \limits_{-x_n < z_n} \delta(\zeta, x, z)K^\alpha(z)dz \\
\leq & -\int \limits_{-x_n < z_n < 0} \delta(\log, x_n, z_n) K^\alpha(z)dz \\
\leq &-\Lambda \int \limits_{-x_n < z_n < 0} \log(1 + z_n/x_n) |z|^{-(n + \alpha)}dz \\
= & -\Lambda \int \limits_{-x_n < z_n < 0} \log(1 + z_n/x_n) |z_n|^{-(n + \alpha)} 
\int \limits_{\R^{n-1}} (1 + |z'/|z_n||^2)^{-(n + \alpha)/2}dz' dz_n.
\end{split}
\end{equation*}

Then, performing the change $y=z'/|z_n|$ and using the integrability of the function $(1 + |y|^2)^{-(n + \alpha)/2}$ in $\R^{n-1}$, we conclude
\begin{equation}\label{changeapp}
-\mathcal{I}_{\Omega}[\zeta](x) \leq  -\Lambda C_{n, \alpha} \int \limits_{-x_n < z_n < 0} \log(1 + z_n/x_n) |z_n|^{-(1 + \alpha)} dz_n.
\end{equation}

Thus, applying the change of variables $t = z_n/x_n$ in the integral of the right-hand side we conclude
\begin{equation*}
-\mathcal{I}_{\Omega}[\zeta](x) \leq -\Lambda C_{n,\alpha} x_n^{-\alpha} \int \limits_{-1}^{0} \log(1 + t)|t|^{-(1 + \alpha)} dt. 
\end{equation*}

Note that in $[-1/2, 0]$ we have
$|\log(1 + t)||t|^{-(1 + \alpha)} \leq 3/4|t|^{-\alpha}$ and since $\alpha < 1$ this term is integrable. Meanwhile in $[-1, -1/2]$ the term 
$|t|^{-(1 + \alpha)}$ is uniformly bounded and by the integrability of the $\log$ function at zero we conclude the lemma for the flat boundary case.

Now we deal with the general case. Take $x \in \Omega$ close to the boundary and consider $\Theta_x$ as in~\eqref{Theta}. By the definition 
of $\zeta$ in~\eqref{zeta} we can use same arguments as at the begining of the proof of Lemma~\ref{lemad}, to conclude
\begin{equation*}
\mathcal{I}_{\Omega}[\zeta](x) \leq \Lambda \int_{\Theta_x} [\log(d(x)) -\log(d(x + z))]|z|^{-(n + \alpha)}. 
\end{equation*}

Denote by $\hat{x}$ the unique point of $\partial \Omega$ such that 
$d(x) = |x - \hat{x}|$. After a rotation we can assume $x - \hat{x} = d(x) e_n$, where $e_n = (0, ..., 0 ,1)$. By the smoothness of the boundary,
there exists $R > 0$, an open set $\mathcal{O} \subset \R^{n-1}$ containing the origin and a smooth function $F: \R^{n - 1} \to \R$ such that 
\begin{equation*}
\{ (x', x_n) \in \mathcal{O} \times \R : x_n = F(x') \} = \partial \Omega \cap B_R(\hat{x}),
\end{equation*}
that is, near $\hat{x}$ the boundary is the graph of a smooth function. We can assume without loss of generality that $F(0) = \hat{x}$. With this, 
we define the change of variables 
\begin{equation*}
\begin{array}{rccl} \Phi : & \mathcal{O} \times [0, d(x)]& \to & \R^n \\ &(y,s) &\mapsto& F(y) - x + (d(x) - s) Dd(F(y)). \end{array}
\end{equation*}

Regarding this function, by the smoothness of $\partial \Omega$, $\Phi$ is a diffeomorphism of $\mathcal{U}$ into $\Phi(\mathcal{U})$,
where $\mathcal{U} \subset \R^n$ is an open set such that $\bar{\mathcal{O}} \times [0,d(x)] \subset \mathcal{U}$. 
This allows us to consider $\Phi$ as a change of variables that flattens $\partial \Omega - x$.
Note that $\Phi$ is $C^1$ and $D\Phi(s,y)$ is uniformly bounded in $\bar{\mathcal{O}} \times [0,d(x)]$,
independent of $x$. Since $x = \hat{x} + d(x)Dd(\hat{x})$, then $\Phi(0) = 0$ and by well-known results, there exists $\kappa > 0$ 
independent of $x$ such that
\begin{eqnarray}\label{estimatekernelapp}
\kappa^{-1} |(y,s)| \leq |\Phi(y,s)| \leq \kappa |(y,s)|,
\end{eqnarray}
for all $(y,s) \in \bar{\mathcal{O}} \times [0, d(x)]$. Denoting
$
\Theta^0_x = \Phi(\mathcal{O} \times [0, d(x)]),
$
we have $\Theta^0_x \subset \Theta_x$ and for each $z \in \Theta_x^0$ there exists a unique 
$(y,s) \in \mathcal{O} \times [0,d(x)]$ such that $z = \Phi(y,s)$ and then $d(x + z) = d(x) - s$. Applying $\Phi$ as a change of variables we can
write 
\begin{equation*}
\begin{split}
& \int_{\Theta_x^0} -\delta(\zeta, x, z)|z|^{-(n+\alpha)}dz \\
&= -\int_{\mathcal{O} \times [0, d(x)]} [\log(d(x) - s) - \log(d(x))]|\Phi(y,s)|^{-(n + \alpha)}|\mbox{det}(D\Phi(y,s))|dyds \\
& \leq -C \kappa^{n + \alpha} \int \limits_{0}^{d(x)} \int \limits_{\mathcal{O}} [\log(d(x) - s) - \log(d(x))] |(y,s)|^{-(n + \alpha)}dyds,
\end{split}
\end{equation*}
where, in the last inequality, we have used the boundedness of $|\mbox{det}(D\Phi)|$ and~\eqref{estimatekernelapp}. Thus, we can integrate over $y$ in a similar
way as in the flat case (see~\eqref{changeapp}) and making $t = -s$ we conclude
\begin{equation*}
\int_{\Theta_x^0} -\delta(\zeta, x, z)|z|^{-(n+\alpha)}dz \leq 
-C \kappa^{n + \alpha} \int \limits_{-d(x)}^{0} \log(1 + t/d(x))|t|^{-(1 + \alpha)}dt, 
\end{equation*}
arriving to the same integral obtained in the flat case. We conclude that the integral over $\Theta_x^0$ is bounded 
by $C d(x)^{-\alpha}$, where $C$ depends on $n, \alpha$ and the smoothness of $\partial \Omega$. 

Since the remaining portion 
$\Theta_x \setminus \Theta_x^0$ is at distance at least $R/2$ from the origin, the term $|z|^{-(n + \alpha)}$ is bounded 
on $\Theta_x \setminus \Theta_x^0$ and then we have
\begin{equation*}
\int_{\Theta_x \setminus \Theta_x^0} -\delta(\zeta, x, z)|z|^{-(n+\alpha)}dz \leq 
C_R \int_{\Theta_x \setminus \Theta_x^0} -\delta(\zeta, x, z)dz,
\end{equation*}
where $C_R$ depends on $R$ but not in $x$. Using similar flattening arguments and the integrability of the $\log$ function near zero
we conclude the last integral is just bounded independently of $x$. This concludes the proof.
\qed

\subsection{Proof of Estimate {\bf \eqref{estimatelemaGammaoutapp}}.} We recall that
\begin{equation*}
\mathcal{A}_{\delta, \mu} = \{ z \in B_r \ : \ d(x) - \delta < d(x + z) < d(x) + \mu\},
\end{equation*}
and that the function $\Psi$ is defined as
$$
\Psi(x) = \eta^{-1} d(x) + \epsilon^{-1}|x - x_0|^2
$$ 
in a neighborhood of $x_0$. Since $\bar{x} \to x_0$, taking $r, \delta$ and $\mu$ suitably small in the definition of $\mathcal{A}_{\delta, \mu}$,
the second integral in~\eqref{estimatelemaGammaoutapp} can be estimated as
\begin{equation}\label{intapp}
\int \limits_{\mathcal{A}_{\delta, \mu}} \delta(\Psi, \bar{x}, z)K^\alpha(z)dz \leq C(\eta^{-1} + \epsilon^{-1}) 
\int \limits_{\mathcal{A}_{\delta, \mu}} |z|^{-(n + \alpha - 1)}dz
\end{equation}
with $C$ independent of $\bar{x}, \eta$ and $\epsilon$. From this point we can use a flattening argument very similar to the previous one, 
building a change of variables which allows us to estimate the whole integral by the one dimensional integral in the normal direction to
the boundary. With this, we can conclude that
\begin{equation*}
\int \limits_{\mathcal{A}_{\delta, \mu}} |z|^{-(n + \alpha - 1)}dz \leq C \int \limits_{-\delta}^{\mu} |s|^{-\alpha}ds,
\end{equation*}
and since we assume $\delta \leq \mu$ then we can bound above the integral by $C\mu^{1 - \alpha}$. This 
concludes~\eqref{estimatelemaGammaoutapp}.
\qed

\medskip

\noindent
{\bf Aknowledgements:} The author would like to express his warmest gratitude to Prof. Guy Barles, for proposing the problem and for its
active participation. Without him, this work would not have been carried out.
The author was partially supported by CONICYT, Grants Capital Humano Avanzado and Realizaci\'on de Tesis Doctoral, and Grant 
Ayudas Para Estad\'ias Cortas of Vicerrectoria de Asuntos Acad\'emicos de la Universidad de Chile.


\begin{thebibliography}{00}

\bibitem{Alvarez-Tourin}
Alvarez, O and Tourin, A. {\em Viscosity Solutions of Nonlinear Integro-Differential Equations} Annales de L'I.H.P., section C, vol.13 (1996),
no. 3, 293-317.


\bibitem{Bardi-Capuzzo}
Bardi, M. and Capuzzo-Dolcetta, I. {\em Optimal Control and Viscosity Solutions of Hamilton-Jacobi-Bellman Equations.} System \& Control:
Foundation and Applications, Birkhauser - Verlag (1997).

\bibitem{Barles-book}
Barles, G. {Solutions de Viscosite des Equations de Hamilton-Jacobi} Collection ``Mathematiques et Applications'' de la SIAM, no 17, Springer-Verlag 
(1994).

\bibitem{Barles-Burdeau}
Barles, G. and Burdeau, J. {\em The Dirichlet Problem for Semilinear Second-Order Degenerate Elliptic Equations and Applications to 
Stochastic Exit Time Control Problems} Comm. in PDE., 20 (1-2), 129-178 (1995). 

\bibitem{Barles-Chasseigne-Imbert}
Barles, G., Chasseigne, E. and Imbert, C. {\em On the Dirichlet Problem for Second Order Elliptic Integro-Differential Equations}
Indiana U. Math. Journal, 2008.

\bibitem{Barles-Chasseigne-Georgelin-Jakobsen}
Barles, G., Chasseigne, E., Georgelin, C. and Jakobsen, E. {\em On Neumann Type Problems for Nonlocal Equations Set in a Half Space.} Preprint.

\bibitem{Barles-DaLio}
Barles, G. and Da Lio, F. {\em Remarks on the Dirichlet and state-constraint problems for quasilinear parabolic equations.}
Advances in Diff. Equations 8 (2003) 897-922.

\bibitem{Barles-Imbert}
Barles, G. and Imbert, C. {\em Second-order Eliptic Integro-Differential Equations: Viscosity Solutions' Theory Revisited.} 
IHP Anal. Non Lin\'eare, Vol. 25 (2008) no. 3, 567-585.

\bibitem{Barles-Rouy}
Barles, G. and Rouy, E. {\em A Strong Comparison Result for the Bellaman Equation Arising in Stochastic Exit Time Control Problems ans its 
Applications.} Comm. in PDE, 23 (11 \& 12)(1998), 1995-2033. 

\bibitem{Caffarelli-Silvestre1}
Caffarelli, L. and Silvestre, L. {\em Regularity Theory For Nonlocal Integro-Differential Equations.} Comm. Pure Appl. Math, Vol. 62 (2009), no. 5, 597-638.


\bibitem{Chasseigne}
Chasseigne, E. {\em The Dirichlet problem for some nonlocal diffusion equations.} Differential Integral Equations 20 (2007), no. 12, 1389–1404. 


\bibitem{usersguide} Crandall, M.G., Ishii H. and Lions, P.-L. {\em User's Guide to Viscosity Solutions of Second Order Partial
Differential Equations.} Bull. Amer. Math. Soc. (N.S.), Vol. 27 (1992), no. 1, 1-67.

\bibitem{DaLio}
Da Lio, F. {\em Comparison Results for Quasilinear Equations in Annular Domains and Applications.} Comm. Partial Diff. Equations, 27 (1 \& 2) 
283-323 (2002).

\bibitem{Hitch}
Di Neza, E., Palatucci, G. and Valdinoci, E. {\em Hitchhiker's Guide to the Fractional Sobolev Spaces.}
 Bull. Sci. Math., 136, (2012), no. 5, 521--573.


\bibitem{Fleming-Soner}
Fleming, W. and Soner, H. {\em Controlled MArkov Processes and Viscosity Solutions} Applications of Mathematics, Springer-Verlag, New York, 1993.

\bibitem{Freidlin}
Freidlin, M.I. {\em Functional Integration and Partial Differential Equations.} Annals of Math. Studies, no 109, Princeton University Press, 1985.

\bibitem{G-T}
D. Gilbarg and N.S. Trudinger, \emph{Elliptic partial differential equations of second order},  Springer-Verlag, Berlin 2001.


\bibitem{Ishii1} Ishii, H. {\em A boundary value problem of the Dirichlet type for Hamilton-Jacobi equations} Ann.
Scuola Norm. Sup. Pisa Cl. Sei. (4) 16 (1989), 105-135.

\bibitem{Ishii2}
Ishii, H. {\em Perron's Method for Hamilton-Jacobi Equations.} Duke Math. J. 55 (1987), 369-384.

\bibitem{Katsoulakis}
Katsoulakis, M. {\em Viscosity Solutions of Second-Order Fully Nonlinear Elliptic Equations with State Constraints.} Indiana Univ. Math. J.,
43, no 2, 493-519, 1994.

\bibitem{Keldysh}
Keldysh, M.V. {\em On some cases of degenerate elliptic equations.} Dokl. Acad. Nauk SSSR, 77, 181-183, 1951.

%
%
%

\bibitem{Radkevich1}
Radkevich, E. V. {\em Equations with Nonnegative Characteristic Form. I.} Journal of Math. Sciences, vol 158, no 3, 297-452, 2009. 

\bibitem{Radkevich2}
Radkevich, E. V. {\em Equations with Nonnegative Characteristic Form. II.} Journal of Math. Sciences, vol 158, no 4, 453-604, 2009. 

\bibitem{Sato}
Sato, K.-I. {\em L\'evy Processes and Infinitely Divisible Distributions.} Cambridge Univ. Press., 1999.

\bibitem{Sayah1}
Sayah, A. {\em \'Equations d'Hamilton-Jacobi du Premier Ordre Avec Termes Int\'egro-Diff\'erentiels. I. Unicit\'e des solutions de viscosit\'e.}
Comm. Partial Differential Equations 16 (1991), 1057-1074.

\bibitem{Sayah2}
Sayah, A. {\em \'Equations d'Hamilton-Jacobi du Premier Ordre Avec Termes Int\'egro-Diff\'erentiels. II. Existence de solutions de viscosit\'e.}
Comm. Partial Differential Equations 16 (1991), 1075-1093.

\end{thebibliography}
\end{document}